\documentclass[10pt]{amsart}
\usepackage{amsbsy}
\usepackage{graphicx,epsfig,subfigure,psfrag}
\begin{document}

\newtheorem{theorem}{Theorem}
\newtheorem{proposition}{Proposition}
\newtheorem{lemma}{Lemma}
\newtheorem{corollary}{Corollary}
\newtheorem{definition}{Definition}
\newtheorem{remark}{Remark}
\newcommand{\tex}{\textstyle}
\numberwithin{equation}{section} \numberwithin{theorem}{section}
\numberwithin{proposition}{section} \numberwithin{lemma}{section}
\numberwithin{corollary}{section}
\numberwithin{definition}{section} \numberwithin{remark}{section}
\newcommand{\ren}{\mathbb{R}^N}
\newcommand{\re}{\mathbb{R}}
\newcommand{\n}{\nabla}
\newcommand{\p}{\partial}
\newcommand{\iy}{\infty}
\newcommand{\pa}{\partial}
\newcommand{\fp}{\noindent}
\newcommand{\ms}{\medskip\vskip-.1cm}
\newcommand{\mpb}{\medskip}
\newcommand{\AAA}{{\bf A}}
\newcommand{\BB}{{\bf B}}
\newcommand{\CC}{{\bf C}}
\newcommand{\DD}{{\bf D}}
\newcommand{\EE}{{\bf E}}
\newcommand{\FF}{{\bf F}}
\newcommand{\GG}{{\bf G}}
\newcommand{\oo}{{\mathbf \omega}}
\newcommand{\Am}{{\bf A}_{2m}}
\newcommand{\CCC}{{\mathbf  C}}
\newcommand{\II}{{\mathrm{Im}}\,}
\newcommand{\RR}{{\mathrm{Re}}\,}
\newcommand{\eee}{{\mathrm  e}}
\newcommand{\LL}{L^2_\rho(\ren)}
\newcommand{\LLL}{L^2_{\rho^*}(\ren)}
\renewcommand{\a}{\alpha}
\renewcommand{\b}{\beta}
\newcommand{\g}{\gamma}
\newcommand{\G}{\Gamma}
\renewcommand{\d}{\delta}
\newcommand{\D}{\Delta}
\newcommand{\e}{\varepsilon}
\newcommand{\var}{\varphi}
\newcommand{\lll}{\l}
\renewcommand{\l}{\lambda}
\renewcommand{\o}{\omega}
\renewcommand{\O}{\Omega}
\newcommand{\s}{\sigma}
\renewcommand{\t}{\tau}
\renewcommand{\th}{\theta}
\newcommand{\z}{\zeta}
\newcommand{\wx}{\widetilde x}
\newcommand{\wt}{\widetilde t}
\newcommand{\noi}{\noindent}
\newcommand{\uu}{{\bf u}}
\newcommand{\xx}{{\bf x}}
\newcommand{\yy}{{\bf y}}
\newcommand{\zz}{{\bf z}}
\newcommand{\aaa}{{\bf a}}
\newcommand{\cc}{{\bf c}}
\newcommand{\jj}{{\bf j}}
\newcommand{\ggg}{{\bf g}}
\newcommand{\UU}{{\bf U}}
\newcommand{\YY}{{\bf Y}}
\newcommand{\HH}{{\bf H}}
\newcommand{\GGG}{{\bf G}}
\newcommand{\VV}{{\bf V}}
\newcommand{\ww}{{\bf w}}
\newcommand{\vv}{{\bf v}}
\newcommand{\hh}{{\bf h}}
\newcommand{\di}{{\rm div}\,}
\newcommand{\ii}{{\rm i}\,}
\newcommand{\inA}{\quad \mbox{in} \quad \ren \times \re_+}
\newcommand{\inB}{\quad \mbox{in} \quad}
\newcommand{\inC}{\quad \mbox{in} \quad \re \times \re_+}
\newcommand{\inD}{\quad \mbox{in} \quad \re}
\newcommand{\forA}{\quad \mbox{for} \quad}
\newcommand{\whereA}{,\quad \mbox{where} \quad}
\newcommand{\asA}{\quad \mbox{as} \quad}
\newcommand{\andA}{\quad \mbox{and} \quad}
\newcommand{\withA}{,\quad \mbox{with} \quad}
\newcommand{\orA}{,\quad \mbox{or} \quad}
\newcommand{\atA}{\quad \mbox{at} \quad}
\newcommand{\onA}{\quad \mbox{on} \quad}
\newcommand{\ef}{\eqref}
\newcommand{\mc}{\mathcal}
\newcommand{\mf}{\mathfrak}

\newcommand{\ssk}{\smallskip}
\newcommand{\LongA}{\quad \Longrightarrow \quad}
\def\com#1{\fbox{\parbox{6in}{\texttt{#1}}}}
\def\N{{\mathbb N}}
\def\A{{\cal A}}
\newcommand{\de}{\,d}
\newcommand{\eps}{\varepsilon}
\newcommand{\be}{\begin{equation}}
\newcommand{\ee}{\end{equation}}
\newcommand{\spt}{{\mbox spt}}
\newcommand{\ind}{{\mbox ind}}
\newcommand{\supp}{{\rm supp}\,}
\newcommand{\dip}{\displaystyle}
\newcommand{\prt}{\partial}
\renewcommand{\theequation}{\thesection.\arabic{equation}}
\renewcommand{\baselinestretch}{1.1}
\newcommand{\Dm}{(-\D)^m}

\title
{\bf The Cauchy problem for a tenth-order thin film equation II.
Oscillatory source-type and \\ fundamental similarity solutions}


\author{P.~\'Alvarez-Caudevilla, J.D.~Evans, and V.A.~Galaktionov}

\address{Universidad Carlos III de Madrid,
Av. Universidad 30, 28911-Legan\'es, Spain -- Work phone number:
+34-916249099} \email{pacaudev@math.uc3m.es}

\address{Department of Mathematical Sciences, University of Bath,
 Bath BA2 7AY, UK -- Work phone number: +44 (0)1225386994 }
\email{masjde@bath.ac.uk}

\address{Department of Mathematical Sciences, University of Bath,
 Bath BA2 7AY, UK -- Work phone number: +44 (0)1225826988}
\email{masvg@bath.ac.uk}


\thanks{This works has been partially supported by the Ministry of Economy and Competitiveness of
Spain under research project MTM2012-33258.}


\date{\today}





\begin{abstract}

 Fundamental global similarity solutions of the standard form
  $$
  \tex{
  u_\g(x,t)=t^{-\a_\g}
 f_\g(y),\,\,\mbox{with the rescaled variable}\,\,\,
 y= x/{t^{\b_\g}}, \,\, \b_\g= \frac {1-n \a_\g}{10},
 }
  $$
where $\a_\g>0$ are real {\em nonlinear eigenvalues}
   ($\g$ is a multiindex in $\ren$)
  of the tenth-order thin
 film equation (TFE-10)
 \begin{equation}
 \label{i1a}
    u_{t} = \nabla \cdot(|u|^{n} \n \D^4 u)
 \quad \mbox{in} \quad \ren \times \re_+
    \,, \quad n>0,
   \ee
  are studied. The present paper continues the study began in \cite{AEG1}. Thus, 
  the following  questions are also under scrutiny:

     {\bf (I)} Further study of the limit $n \to 0$, where the
     behaviour of finite interfaces and solutions  as $y \to \iy$ are
     described. In particular, for $N=1$, the interfaces are shown to diverge
     as follows:
     $$
      \tex{
     |x_0(t)| \sim 10 \left( \frac{1}{n}\sec\left( \frac{4\pi}{9} \right) \right)^{\frac 9{10}}  t^{\frac 1{10}}
      \to \iy \asA n \to 0^+.
     }
      $$

      {\bf (II)} For a fixed $n \in (0, \frac 98)$, oscillatory structures of solutions
      near
      interfaces.

      {\bf (III)} Again, for a fixed $n \in (0, \frac 98)$, global structures of some nonlinear
      eigenfunctions $\{f_\g\}_{|\g| \ge 0}$ by a combination of
      numerical and analytical methods.

\end{abstract}

\maketitle

\section{The TFE-10 and  a nonlinear eigenvalue problem }
 \label{S1}

\subsection{Main model and previous  results}

\noindent
    We study the global-in-time behaviour of compactly supported solutions of the Cauchy problem of
a tenth-order quasilinear evolution equation of parabolic type,
called the {\em thin film equation} (TFE--10)
\begin{equation}
\label{i1}
    u_{t} = \nabla \cdot(|u|^{n} \n \D^4 u)
 \quad \mbox{in} \quad \ren \times \re_+
    \,,
\end{equation}
where $\n={\rm grad}_x$ and $n>0$ is a real parameter. In view of
the degenerate mobility coefficient $|u|^n$, equation \eqref{i1}
is written for solutions of changing sign, which can occur in the
{\em Cauchy problem} (CP) and also in some {\em free boundary
problems} (FBPs).

Equation \ef{i1} has been chosen as a typical difficult
higher-order quasilinear degenerate parabolic model.
Though the fourth-order version  has been the most studied
actually, since 1980s (and, later on, in the 2000s, various
six-order ones), higher-order
quasilinear degenerate equations are known to occur in several
applications and, during the last ten-fifteen years, have began to
steadily penetrate into modern nonlinear PDE theory; see a number
of references/results in \cite[\S~1.1]{EGK3} and in
\cite{CG,Liu,LQ}. Concerning the origin of the TFE--10 \ef{i1}, as
in \cite{AEG1}, honestly, we have chosen this very difficult model
in order to develop mathematical PDE techniques showing that such
complicated quasilinear degenerate  equations, anyway, admit a
rather constructive study {\em regardless} their order. As is
well-known, the spatial order of nonlinear PDEs essentially affect
their difficulty for any proper mathematical study. However, our
techniques are now not {\em order-dependent}, so the many results
can be extended to any TFE--12, --14,\,...\,,--20, etc. Of course,
in such a generality of the approaches, we do not pretend to give a
full description of all the solution properties, though notice
several common properties of all types of TFEs.

\ssk

For convenience, first, we need to  state the main result obtained
in our previous paper \cite{AEG1}, a study to be continued here.
 Thus, in \cite{AEG1},
we introduced global self-similar solutions of \ef{i1} of the
standard form
\begin{equation}
\label{sf3}
 \tex{
    u(x,t):= t^{-\a} f(y), \quad \hbox{with}\quad
    y:=\frac{x}{t^\b},\quad \b= \frac{1-n \a}{10},
    }
\end{equation}
 where $\a>0$ stands for the so-called real {\em nonlinear eigenvalues} and the {\em nonlinear eigenfunctions}  $f$ satisfy an elliptic equation
 \begin{equation}
\label{self1}
 \tex{
     \nabla \cdot \left( |f|^{n} \n \D^4 f\right) +\frac{1-\a n}{10}\, y \cdot \nabla f +\a
    f=0,
    \quad f \in C_0(\ren)\, .
    }
\end{equation}

 Then, we state a  {\em nonlinear eigenvalue problem}
  for pairs $\{\a,f\}$\footnote{More precisely, as in \cite{AEG1}, we actually deal with
a  ``linear (in $\a$) spectral pencil for the quasilinear TFE-10
 operator".},
 where the problem setting includes finite propagation phenomena
 for such TFEs, i.e., $f$ is assumed to be compactly supported, $f \in
 C_0(\ren)$. This is a kind of an assumed ``minimal"
behaviour of $f(y)$ as $y \to \iy$, which naturally accompany many
standard singular Sturm--Liouville problems and others.

Using long-established terminology, we call such similarity
solutions \ef{sf3} (and also the corresponding profiles $f$) to be
a sequence of  {\em fundamental solutions}. Though, actually, the
classic fundamental solution is the first radially symmetric one
(with the first kernel $f_0=f_0(|y|)$), which is  the {\em
instantaneous source-type solution} of \ef{i1} with Dirac's delta
as initial data. Moreover, for $n=0$, $f_0(|y|)$ becomes the
actual rescaled kernel of the fundamental solution  of the linear
operator $D_t-\D_x^5$.

\ssk

Our main goal in \cite{AEG1} was  to show {\em analytically} that,
at least, for small $n>0$,
 \be
 \label{main1}
 \mbox{(\ref{self1}) admits a countable set of fundamental
 solutions $\Phi(n)=\{\a_\g,f_\g\}_{|\g| \ge 0}$},
  \ee
  where $\g$ is a multiindex in $\ren$ to numerate
   these eigenvalue-eigenfunction pairs.

\ssk

 Studying the nonlinear eigenvalue problem \ef{self1} in
\cite{AEG1}, we  performed a ``homotopic deformation" of \ef{i1}
as $n \to 0^+$
and reduced it to the classic {\em poly-harmonic
equation of tenth order}
\begin{equation}
\label{lin5}
    \tex{ u_{t} = \D^5 u\quad \hbox{in} \quad \re^{N} \times \re_+\,.}
\end{equation}

\ssk

 More precisely, we answered in \cite{AEG1} the following
question:

\ssk

 {\bf (I)}  Passing
    to the limit $n \to 0^+$ in \ef{i1a} on any compact subsets $\{|y| \le C\}$
      by using Hermitian non-self-adjoint
    spectral theory for a pair of rescaled non-symmetric operators $\{\BB,\BB^*\}$ corresponding to the  linear
    poly-harmonic equation \ef{lin5} where
     \be
     \label{Boper}
      \tex{
     \BB=\D^5 + \frac 1{10}\, y \cdot \n+
     \frac N{10}I \andA \BB^*=\D^5 - \frac 1{10}\, y \cdot \n\, .
     }
     \ee
The corresponding problem \ef{self1} then reduces to a standard
(but not self-adjoint) Hermitian type  linear eigenvalue problem
for the pair $\{\BB,\BB^*\}$. Therefore, according to this
approach, the nonlinear version of
 \ef{main1} has the origin in the discreteness-reality of the
 spectrum of the corresponding linear operator $\BB$.

     This allowed us to identify a countable family of nonlinear
     eigenfunctions for \ef{i1a}, at least, for small $n>0$, which
     defined proper solutions of the Cauchy problem for the
     TFE-10.

\subsection{Main new results and layout of the paper}

 In the present paper, our main ``non-local" goal
    is to verify a possibility of global extensions
  of such ``$n$-branches" of some first fundamental solutions, which was   then checked
  numerically. A couple of such  preliminary results were already available in \cite[\S~5]{AEG1}.

  To do so, first we analyse the limiting behaviour of the problem \eqref{i1} in the one dimensional case, obtaining
  an approximating structure for the solutions satisfying \eqref{i1}
   close to the interfaces. Also, in one-dimension, we ascertain, via numerical and analytical
   methods,
  the existence of periodic oscillatory structure of solutions for $n\in (0,\frac{9}{8})$.

  Furthermore, through a homotopic approach, and using
   the standard degree theory, we ascertain the existence of a countable family of global similarity profiles for
   \eqref{i1}.

Section \ref{S2} is devoted to similarity solutions and derivation
of the corresponding nonlinear eigenvalue problem. Later, we
address the following questions:

\ssk

     {\bf (II)} Section \ref{S3}: Further study of the limit $n \to 0$, where the
     behaviour of finite interfaces and solutions  as $y \to \iy$ are
     described. Analysing the limiting behaviour when $n\to 0^+$ in the one dimensional case we obtain that
     the non-uniform solution in this limit comprises two regions
     $$
     \hbox{an \rm{Inner region}}\quad \{x=O(1)\}\quad
     \hbox{and}\quad \hbox{an \rm{Outer region}}\quad \{x=O(n^{-9/10})\},
     $$
in which $u$ is exponentially small. In the \rm{Inner region}
$|u|^n \sim 1$ for small $n>0$. Hence, the asymptotic behaviour of
the solution tells us that the solution satisfies the
poly-harmonic equation
   $$
   \tex{
    \frac{\partial u_0}{\partial t} = \frac{\partial^{10} u_0}{\partial x^{10}}
    \, .
    }
    $$
    Moreover, from the performed analysis, it is clear that this solution breaks down when $\{x=O(n^{-9/10})\}$, the \rm{Outer region}.
     In particular, the interfaces $x=x_0(t)$ are shown to diverge
     as follows:
     $$
     \tex{
     |x_0(t)| \sim 10 \left( \frac{1}{n}\sec\left( \frac{4\pi}{9} \right) \right)^{\frac 9{10}}  t^{\frac 1{10}}  \to \iy \asA n \to 0.
     }
      $$
      By a similar analysis, we also obtain the structure of the eigenfunctions satisfying
      \be
 \label{ODE1}
  \tex{
 (|f|^n f^{(9)})' +\frac{1-\a n}{10}\, y \, f' +\a
    f=0,
    \quad f \in C_0(\re)\, .
    }
    \ee

\ssk

      {\bf (III)} Section \ref{S4}: For a fixed $n \in (0,\frac 98)$, oscillatory structures that occur near
      interfaces are detected. In one dimension, we study the local behaviour near the interface $y=\frac{x_0}{t^\b}$ for the similarity ODE \eqref{ODE1}
    assuming similarity profiles $f(y)$ with a \rm{maximal regularity} at the interface $y=y_0$, allowing the trivial extension
     $f=0$ for $y  \ge y_0$.
        An analysis of periodic solutions for small $n>0$ is also performed.
\ssk

      {\bf (IV)} Section \ref{S5}: Again, for a fixed $n \in (0, \frac 98)$, the study of global structures of some nonlinear
      eigenfunctions $\{f_\g\}_{|\g| \ge 0}$ by a  combination of
      numerical and analytical methods. Using a homotopy transformation of the form
      $$
      |f|^n \mapsto |\e^2+f^2|^{\frac{n}{2}}, \quad \e>0\,\,\,
      \mbox{small},
      $$
      and applying the standard degree theory, we perform a
            double limit when $\e , n\to 0^+$ obtaining existence and multiplicity results for the oscillatory solutions of changing
      sign of the non-linear elliptic equation \eqref{self1}. Then,
            we arrive at the existence of a countable family of solutions by a direct $n$-expansion from
      the solution of the linear elliptic equation whose operator is denoted by \eqref{Boper}.

Finally, in Section \ref{S6}, we construct some of those nonlinear
eigenfunctions for the non-linear
            elliptic equation \eqref{self1} numerically in the one-dimensional case. This analysis provides us with a graphic idea of
            the behaviour of the solutions obtained in the previous sections.

\subsection{
Possible origins of discrete
nonlinear spectra and principle difficulties}

It is key for us that \ef{self1} {\em is not variational}, so we
cannot use powerful tools such as Lusternik--Schnirel'man (L--S,
for short) category-genus theory, fibering, and other well-known
approaches, which in many cases are known to provide at least a
{\em countable} family of critical points (solutions) in the
non-coercive case, when the category of the functional subset
involved is typically infinite.

It is also crucial and well known that the L--S min-max approach
{\em does not detect all families of critical points}. However,
sometimes it can revive some amount of solutions. A somehow
special example was revealed   in \cite{GMPSobI, GMPSobII}, where
key features of those variational L--S and fibering approaches
applied are described. Namely,  for some variational fourth-order
and higher-order ODEs in $\re$, including  those with the typical
non-linearity $|f|^n f$, as above,
  \be
 \label{mm.561}
  \mbox{$
 -  (|f|^n f)^{(4)} + |f|^n f= \frac 1n \, f \inB  \re, \quad f \in C_0(\re) \quad (n>0),
 $}
  \ee
 as well as for the following standard looking
one with the only cubic nonlinearity \cite[\S~6]{GMPSobII}:
 \be
 \label{anal1}
-f^{(4)}+ f = f^3 \quad \mbox{in} \quad \re, \quad f \in
H^4_\rho(\re) \quad (\rho={\mathrm e}^{a|y|^{4/3}}, \,\, a>0
\,\,\, \mbox{small}).
 \ee
  It was shown then that these equations admit a {\em countable set of
countable families of solutions}, while the L--S/fibering approach
detects only {\sc one} such a discrete family of (min-max)
critical points. Further countable families are not expected to be
determined easily by more advanced techniques of potential theory,
such as the mountain pass theorem,  fibering methods,  or others.
Existence of other, not L--S type critical points for \ef{mm.561}
and \ef{anal1}, were shown in \cite{GMPSobI, GMPSobII} by using a
combination of numerical and (often, formal) analytic methods and
heavy use of oscillatory nature of solutions close to finite
interfaces (for \ef{mm.561}) and at infinity (for \ef{anal1}). In
particular, detecting the corresponding L--S countable sequence of
critical points was done {\em numerically}, i.e., by checking
their actual min-max features (their critical values must be
maximal among other solutions belonging to the functional subset
of a given category, and having  a ``suitable geometric shape").

 Therefore, even in
the variational setting, counting various families of critical
points and values represents a difficult open problem for such
higher-order ODEs, to say nothing of their elliptic counterparts
in $\ren$.

Hence,
in \cite{AEG1}, we  relied on a different approach, in particular,
a ``homotopic deformation" of \ef{i1}  as $n \to 0^+$, which is
also effective for such difficult variational problems and detects
more solutions than L--S/fibering theory (though only locally upon
the parameter). The philosophy  in which that ``homotopic
approach" was performed in \cite{AEG1}, is the base of our
analytic approach of section\,\ref{S5} in order to get that
countable family of solutions.

\subsection{The second model: bifurcations in $\re^2$}

 Our homotopy approach can be extended to a
 more complicated   {\em unstable thin film equation} (TFE--10) in
 the critical case
\begin{equation}
\label{e1}
 \tex{
    u_{t} = \nabla \cdot(|u|^{n} \n \D^4 u)-\D(|u|^{p-1}u)
 \quad \mbox{in} \quad \ren \times \re_+
    \,, \quad p>n+1,
    }
\end{equation}
with the extra unstable diffusion term. We obtain
a discrete real nonlinear spectrum for
\ef{e1} that requires a simultaneous {\em double} homotopy
deformation $n \to 0^+$ and $p \to 1^+$ leading to a new linear
Hermitian spectral theory. We do not develop it here and just
focus on a principal opportunity to detect a discrete nonlinear
spectrum for \ef{e1}. More details on blow-up and global
similarity solutions (as unique extensions after blow-up) of
\ef{e1} can be found in \cite{AEG3}.


\subsection{Global extension of bifurcation branches: a principal
open problem}

It is worth mentioning that, for both problems \ef{self1} and the
corresponding problem occurring for \ef{e1} (after the similarity
time-scaling), a global extension of bifurcation
$n$-branches ($(n,p)$-branches for \ef{e1}) represents a difficult
open problem of general nonlinear operator theory. Moreover, as
was shown in \cite{GalPetII} (see also other examples in
\cite{GMPSobII}), the TFE-4 with absorption $-|u|^{p-1}u$ (instead
of the backward-in-time diffusion as in \ef{self1}), depending on
not that small $n \sim 1$, admits some $p$-bifurcation branches
having turning (saddle-node) points and thus representing  closed
loops. Hence, these branches {\em are not globally extendable} in
principle. On the other hand, for equations with monotone
operators such as the PME-4
\begin{equation}
 \label{PME4}
    u_t =-(|u|^{n} u)_{xxxx}\quad \hbox{in} \quad \re \times \re_+\,,
\end{equation}
the
$n$-branches seem to be globally extensible in $n>0$,
\cite{GalRDE4n}.

\setcounter{equation}{0}
\section{Problem setting and self-similar solutions}
\label{S2}


\subsection{The FBP and CP}

As done previously in \cite{EGK1}--\cite{EGK4}, we distinguish the standard
{free-boundary problem} (FBP) for \ef{i1} and the {\em Cauchy
problem}; see further details therein.

 For  both the FBP and the
CP, the solutions are assumed to satisfy standard free-boundary
conditions or boundary conditions at infinity:
\begin{equation}
\label{i3}
    \left\{\begin{array}{ll} u=0, & \hbox{zero-height,}\\
    \nabla u=\nabla^2 u=\nabla^3 u=\nabla^4 u=0, & \mbox{``zero-angle",}   \\
    -{\bf n} \cdot (|u|^{n} \nabla \Delta^4 u)=0, &
    \hbox{conservation of mass (zero-flux)}\end{array} \right.
\end{equation}
at the singularity surface (interface) $\Gamma_0[u]$, which is the
lateral boundary of
\begin{equation}
 \label{gamma1}
    \hbox{supp} \;u \subset \ren \times \re_+,\quad N \geq 1\,,
\end{equation}
where ${\bf n}$ stands for the unit outward normal to
$\Gamma_0[u]$.  Note that, for sufficiently smooth interfaces,
the condition on the flux can be read as
\begin{equation*}
    \lim_{\hbox{dist}(x,\Gamma_0[u])\downarrow 0}
    -{\bf n} \cdot \nabla (|u|^{n}  \Delta^4 u)=0.
\end{equation*}
This condition is directly related with the conservation of mass.

Moreover, we also assume bounded, smooth, integrable,
compactly supported initial data
\begin{equation}
\label{i4}
    u(x,0)=u^0(x) \quad \hbox{in} \quad \Gamma_0[u] \cap \{t=0\}.
\end{equation}

\ssk

For the CP, the assumption of nonnegativity is got rid of, and
solutions become oscillatory close to interfaces. It is then key,
for the CP, that the solutions are expected to be
``smoother" at the interface than those for the FBP, i.e., \ef{i3}
are not sufficient to define their regularity. These {\em maximal
regularity} issues for the CP, leading to oscillatory solutions,
are under scrutiny in \cite{EGK2} for a fourth-order case.

In the CP for \eqref{i1} in $\ren \times \re_+$, one needs to pose
bounded compactly supported initial data \eqref{i4} prescribed in
$\ren$.

\subsection{Global similarity solutions:  a nonlinear eigenvalue problem}

 We now specify the self-similar
solutions of the equation \eqref{i1}, which
are admitted due to its natural scaling-invariant  nature.
In the case of the mass being conserved, we have global in time
source-type solutions
$$\tex{
    u(x,t):= t^{-\alpha} f(\frac{x}{t^\b}), \quad \alpha= \frac{1-10\b}{n},
    }
    $$
with $f$ solving the quasilinear elliptic equation (nonlinear eigenvalue problem)
given in \ef{self1}. We add to the elliptic equation  a natural assumption that
 $f$ must be compactly supported (and, of course, sufficiently
 smooth at the interface, which is an accompanying question to be
 discussed as well).
For further details of how to obtain them see \cite{AEG1}.

Thus, for such degenerate elliptic equations,
  the functional setting of \ef{self1} assumes that we are
 looking for  (weak) {\em compactly supported} solutions $f(y)$ as
 certain ``nonlinear eigenfunctions" that hopefully occur for special values of nonlinear eigenvalues
  $\{\a_\g\}_{|\g| \ge 0}$. Therefore,  our goal is to justify
  that
 \ef{main1} holds.

 Concerning  the well-known properties of finite propagation for TFEs, we refer to papers
 \cite{EGK1}--\cite{EGK4}, where a large amount of earlier
 references are available; see also \cite{GMPSobI, GMPSobII} for more recent
 results and references in this elliptic area.

 However, one should observe that there are still
a few entirely rigorous results, especially those that are
attributed to the Cauchy problem for TFEs.

In the linear case $n=0$,
 the condition $f \in C_0(\ren)$, is naturally replaced by the requirement that the
 eigenfunctions $\psi_\b(y)$ exhibit typical exponential decay at
 infinity, a property that is reinforced by introducing  appropriate weighted $L^2$-spaces. Complete details about the spectral theory
 for this linear problem when $n=0$ in \cite{EGKP}.
Actually,
 using the homotopy limit $n \to 0^+$, we will be obliged for
 small $n>0$,
 instead of $C_0$-setting in
\eqref{self1}, to use the following weighted $L^2$-space:
  \begin{equation}
   \label{WW11}
  f \in L^2_\rho(\ren), \quad \mbox{where} \quad \rho(y)={\mathrm e}^{a |y|^{10/9}},
  \quad a>0 \,\,\,\mbox{small}.
  \end{equation}

Note that, in the case of the Cauchy problem with conservation of
mass making use of the self-similar solutions \eqref{sf3}, and performing similar computations as done in \cite{AEG1} we have
that
\begin{equation}
\label{alb1}
 \tex{
 -\a + \b N=0 \LongA
    \a_0(n)=\frac{N}{10+Nn} \quad \mbox{and}  \quad \b_0(n)=\frac{1}{10+Nn}.
    }
\end{equation}

\setcounter{equation}{0}
\section{The limit $n \to 0$: Behaviour of finite interfaces and
nearby solutions}
\label{S3}

We consider here the singular limit $n\to 0^+$ for the full
equation (\ref{i1}) in one space dimension $N=1$. The non-uniform
solution in this limit comprises two regions, an {\em Inner
region} $\{x=O(1)\}$, where $u=O(1)$, and an {\em Outer region}
$\{x=O(n^{-\frac 9{10}})\}$, in which $u$ is exponentially small.
The labelling of these regions as inner and outer becomes clearer
during the course of the scalings.

We begin with the region $\{x=O(1)\}$, for which $u=O(1)$ and
consequently $|u|^n \sim 1$ for small $n>0$. At leading order
$u\sim u_0(x,t)$ satisfies the linear poly-harmonic equation
\begin{equation}
 \tex{
    \frac{\partial u_0}{\partial t} = \frac{\partial^{10} u_0}{\partial x^{10}} ,
    }
\label{s3eq1}
\end{equation}
where $u_0$ here is the leading order term in an expansion with respect to $n$
 (and is not the initial function  in (\ref{i4})). We are interested in an oscillatory class of solutions that are analytic in
 $x$. The far-field behaviour of (\ref{s3eq1}) may be determined using a WKBJ expansion in the form
\begin{equation}
   u_0 \sim a(x,t) {\mathrm e}^{-\phi(x,t)} \hspace{0.6cm} \mbox{as \,\,$x \to +\infty$},
\label{eq:s3eq2}
\end{equation}
which gives
\begin{equation}
 \tex{
   \frac{\partial \phi}{\partial t} =  \left( \frac{\partial \phi}{\partial x} \right)^{10} \andA
     \frac{\partial a }{\partial t} + 10 \left( \frac{\partial \phi}{\partial x} \right)^9 \frac{\partial a}{\partial x}   = -45  \left( \frac{\partial \phi}{\partial x} \right)^8 \frac{\partial^2 \phi}{\partial x^2} a . \label{s3eq3}
     }
\end{equation}
The required solutions to (\ref{s3eq3}) take the form
\begin{equation}
 \tex{
       \phi(x,t) = \phi_{\pm}(x,t) \equiv \frac{9}{10^{10/9}}
        {\mathrm e}^{\pm \frac{4\pi i}{9}} \frac{x^{10/9}}{t^{1/9}}, \hspace{.6cm}
             a = t^{-\frac 12} \Psi_{\pm}\big(\frac xt\big),
             }
\label{s5eq4}
\end{equation}
where $\Psi_{\pm}(\zeta)$ are arbitrary smooth functions (depending on the initial data), but satisfy
 $\Psi_{-}(\zeta)=\bar{\Psi}_{+}(\zeta)$. Thus,
\begin{eqnarray}
   && { }\hspace*{-1cm}
    \tex{
    u_0 \sim \frac{1}{t^{1/2}} \Psi_{+} \big(\frac{x}{t} \big)  \exp \big\{- \frac{9}{10^{10/9}}
    {\mathrm e}^{\frac{4\pi i}{9}} \big(\frac{x^{10}}{t} \big)^{\frac{1}{9}}
    \big\}
    }
     \label{s3eq5} \\
    && \hskip 1cm +
      \tex{
      \frac{1}{t^{1/2}} \Psi_{-} \big(\frac{x}{t} \big)  \exp \Big(- \frac{9}{10^{10/9}} {\mathrm e}^{-\frac{4\pi i}{9}} \big(\frac{x^{10}}{t} \big)^{\frac{1}{9}} \Big) \hspace{0.5cm}
    \mbox{as \,\, $x \to +\infty$}.
    }
     \nonumber
\end{eqnarray}
It is clear from (\ref{s3eq5}) that this solution breaks down when
$x=O(n^{-9/10}) $, since we can no longer approximate $|u|^{n}$ by
unity. This suggests the consideration of an outer region with
scaling $X=n^{9/10}x$. In $X=O(1)$, the PDE becomes
\begin{equation}
 \tex{
    \frac{\partial u}{\partial t} =
     n^{9} \frac{\partial}{\partial X} \big( |u|^n \frac{\partial^9 u}{\partial X^9}\big)  ,
    }
\label{s3eq6}
\end{equation}
this being a conventional formulation of a singular problem, where
the small parameter multiplies the highest derivative. However, as
for the fourth- and sixth-order cases (see \cite{EGK2,EGK4}),
there are fast oscillations superposed on the slow exponential
decay that occurs over this length scale, necessitating the
application of a multiple scales (Kuzmak) approach. As such we
introduce the fast variable
 $$
 \tex{
 Z=\frac{\sigma(X,t)}n,
 }
  $$
    where $\sigma(X,t)$ will be
determined in the standard way by the criterion that the
dependence on $Z$ is periodic of constant (rather than
$(X,t)$-dependent) periodicity - without loss of generality, we
take the period to be $2\pi$. The multiple-scales ansatz for this
region takes the form
\begin{equation}
      u \sim {\mathrm e}^{-\Phi(X,t)/n} A(X,Z,t) \hspace{0.6cm} \mbox{as\,\, $n \to 0$},
\label{s3eq7}
\end{equation}
to within an algebraic power of $n$ (which is determined by the far-field behaviour of $\Psi_{\pm}(\zeta)$), wherein
$\Phi$ is real. Thus, as $n \to 0$,
\begin{equation}
 \tex{
     \frac{\partial u}{\partial t} \sim \frac{1}{n} \left( - \frac{\partial \Phi}{\partial t} A + \frac{\partial \sigma}{\partial t}
\frac{\partial A}{\partial Z} \right) {\mathrm e}^{-\Phi / n} ,
}
\label{s3eq8}
\end{equation}
\begin{equation}
 \tex{
 n^9 \frac{\partial}{\partial X} \big( |u|^{n} \frac{\partial^9 u}{\partial X^9}   \big)
 \sim \frac{1}{n} {\mathrm e}^{-\Phi}
  \Big[ \sum_{k=0}^{10}
 \binom {10}k
   \frac{\partial^k A}{\partial X^k} \big( \frac{\partial \sigma}{\partial X} \big)^k
    \big( - \frac{\partial \Phi}{\partial X} \big)^{10-k} \Big] {\mathrm e}^{-\Phi / n}.
  }
\label{s3eq9}
\end{equation}
We remark that these expansions need to be taken to next (i.e.,
$O(n)$ smaller) order if we are to characterise the dependence of
$A$ on $X$ and $t$; we shall not proceed with such an analysis
here. Viewing the balance (\ref{s3eq8}) and (\ref{s3eq9}) as an
ordinary differential equation in $Z$, we observe that the
condition of $2\pi$ periodicity in $Z$ requires that
\begin{equation}
       A = \alpha_{+}(X,t) {\mathrm e}^{iZ} + \alpha_{-}(X,t) {\mathrm e}^{-iZ}
\label{s3eq10}
\end{equation}
with $\alpha_{-}=\bar{\alpha}_{+}$. Grouping real and imaginary parts, we obtain a coupled system for $\Phi$ and $\sigma$ given by the equations
\begin{eqnarray}
&&
 \tex{
 \frac{\partial \Phi}{\partial t} = {\mathrm e}^{- \Phi} \left[  \left( \frac{\partial \sigma}{\partial X} \right)^{10}
   - 45 \left( \frac{\partial \sigma}{\partial X} \right)^8  \left( \frac{\partial \Phi}{\partial X} \right)^2
   + 210 \left( \frac{\partial \sigma}{\partial X} \right)^6  \left( \frac{\partial \Phi}{\partial X} \right)^4 \right.
   }
    \label{s3eq11} \\
&&  \hspace*{2cm}
 \tex{
 \left.
- 210 \left( \frac{\partial \sigma}{\partial X} \right)^4  \left( \frac{\partial \Phi}{\partial X} \right)^6
   + 45 \left( \frac{\partial \sigma}{\partial X} \right)^2  \left( \frac{\partial \Phi}{\partial X} \right)^8
   - \left(\frac{\partial \Phi}{\partial X} \right)^{10} \right] ,
   }
    \nonumber
\end{eqnarray}
\begin{eqnarray}
 &&
  \tex{
   \frac{\partial \sigma}{\partial t} =  {\mathrm e}^{- \Phi} \left[
   - 10 \left( \frac{\partial \sigma}{\partial X} \right)^9 \frac{\partial \Phi}{\partial X}
   + 120 \left( \frac{\partial \sigma}{\partial X} \right)^7  \left(\frac{\partial \Phi}{\partial X} \right)^3
   - 252 \left( \frac{\partial \sigma}{\partial X} \right)^5 \left( \frac{\partial \Phi}{\partial X} \right)^5 \right.
   }
   \label{s3eq12} \\
&& \hspace*{2cm}
\tex{
\left.  + 120 \left( \frac{\partial \sigma}{\partial X} \right)^3 \left(\frac{\partial \Phi}{\partial X} \right)^7
   - 10  \frac{\partial \sigma}{\partial X}  \left(\frac{\partial \Phi}{\partial X} \right)^9
    \right] .
    }
     \nonumber
\end{eqnarray}
Matching to (\ref{s3eq5}) suggests seeking a consistency relation between (\ref{s3eq11}) and (\ref{s3eq12}) of the form $\sigma = \lambda \Phi$ with
$\lambda$ real, leading to
\[
 \lambda^{10}-35\lambda^8 + 90 \lambda^6 + 42\lambda^4-75\lambda^2+9=0\,\,\Longrightarrow \,\,
  (\lambda^2+1)(\lambda^2-3)(\lambda^6 -33
  \lambda^4 + 27\lambda^2 -3)=0 .
\]
The appropriate root of this characteristic equation is
\begin{equation}
 \tex{
\lambda = \tan \big( \frac{4 \pi}{9}\big),
}
\label{s3lambda}
\end{equation}
this being consistent with the ratio of the imaginary to real
parts in the  exponentials in (\ref{s3eq5}). Consequently, we
obtain a Hamilton--Jacobi equation of the form
\begin{equation}
 \tex{
  \frac{\partial \Phi}{\partial t}
          = - \sec^9\left( \frac{4\pi}{9}\right) {\mathrm e}^{- \Phi}  \left( \frac{\partial \Phi}{\partial X} \right)^{10}  ,
          }
\label{s3eq13}
\end{equation}
the required solution being
\[
 \tex{
    \Phi(x,t) = - 9 \ln \big( 1 - \cos\big( \frac{4\pi}{9} \big)  \big( \frac{X^{10}}{10^{10} t} \big)^{\frac{1}{9}}  \big) ,
    }
\]
which matches successfully with the real part of (\ref{s3eq5})  in
the limit $X \to 0$. Thus, the leading order solution in this
region takes the form
\begin{equation}
 \tex{
     u(x,t) \sim \Big( 1 - \cos\big( \frac{4\pi}{9} \big) \big(\frac{X^{10}}{10^{10 } t} \big)^{\frac{1}{9}}   \Big)^{\frac 9n} A(X,Z,t) ,
     }
\label{s3eq14}
\end{equation}
with $A$ as given in (\ref{s3eq10}), this local behaviour having
the expected $\frac 9n$ power-law form with oscillations
superimposed as in (\ref{s3eq7}) and (\ref{s3eq10}). The interface
$x=x_0(t)$ is thus given by
\begin{equation}
 \tex{
     x_0(t) \sim 10 \Big( \frac{1}{n}\sec\big( \frac{4\pi}{9} \big) \Big)^{\frac 9{10}}  t^{\frac 1{10}} \hspace{0.6cm}
     \mbox{ as \,\,$n \to 0$ } ,
     }
\label{s3eq15}
\end{equation}
illustrating its behaviour for small $n$.

We may also determine the structure of the eigenfunctions
satisfying (\ref{ODE1}) in the small $n$ limit. Again, we have a
two region structure: an inner region $\{y=O(1)\}$, in which
$f=O(1)$ and an outer region $\{y=O(n^{-9/10})\}$ where $f$ is
exponentially small. In the inner region $\{y=O(1)\}$, we obtain
at leading order ($f\sim f_0$) in $n$ the linear ODE
\begin{equation}
 \tex{
f_0^{(10)} + \frac{1}{4}\,yf_0' + \alpha f_0 = 0. } \label{s3eq16}
\end{equation}
An explicit general solution can be expressed in terms of
hypergeometric functions, easily obtained using e.g. {\tt Maple}.
The far-field behaviour of (\ref{s3eq1}) may be determined using a
WKBJ expansion in the form
\begin{equation}
   f_0(y) \sim a(y) {\mathrm e}^{-\phi(y)} \hspace{0.6cm} \mbox{as\,\, $y \to +\infty$}
\label{eq:s3eq17}
\end{equation}
which gives
\begin{equation}
    10 (\phi')^9 = y, \hspace{.6cm}
     \left( \alpha - 45 (\phi')^8 \phi''\right) a + \left( \frac y{10} - 10 (\phi')^9 \right) a' = 0 . \label{s3eq18}
\end{equation}
The required solutions to (\ref{s3eq18}) take the form
\begin{equation}
 \tex{
       \phi(y) = \phi_{\pm}(y) \equiv  \frac{9}{10^{10/9}} {\mathrm e}^{\pm \frac{4\pi i}{9}} y^{10/9} , \hspace{.6cm}
             a(y) = k_{\pm} y^{5(2\alpha-1)/9} ,
             }
\label{s5eq19}
\end{equation}
with $k_{\pm}$ arbitrary constants. Thus,
\begin{eqnarray}
   && { }\hspace*{-1cm}
   \tex{
   f_0(y) \sim k_{+} y^{5(2\alpha-1)/9}  \exp \big\{- \frac{9}{10^{10/9}}
    {\mathrm e}^{\frac{4\pi i}{9}} y^{\frac{10}{9}} \big\}
   }\label{s3eq20} \\
    && \hskip 1cm
    \tex{
    +  k_{-} y^{5(2\alpha-1)/9}  \exp \big\{- \frac{9}{10^{10/9}} {\mathrm e}^{-\frac{4\pi i}{9}}
    y^{\frac{10}{9}} \big\} \hspace{0.5cm}
    \mbox{as \,\,$y \to +\infty$}.
    }
     \nonumber
\end{eqnarray}
Again, this solution breaks down when $y=O(n^{-9/10}) $,
suggesting the consideration  of an outer region with scaling
$Y=n^{9/10}y$. In $Y=O(1)$, we have
\begin{equation}
\tex{
    n^{9} \frac{\mathrm d}{{\mathrm d}Y} \big( |f|^n
    \frac{{\mathrm d}^9 f}{{\mathrm d}Y^9} \big) + \big( \frac{1-\alpha n}{10} \big) y \frac{{\mathrm d}f}{{\mathrm d}Y}
     + \alpha f = 0 .
    }
\label{s3eq21}
\end{equation}
Rather than posing a multiple-scales ansatz directly, we may instead consider
\begin{equation}
      f(Y) \sim {\mathrm e}^{b(Y)/n} B(Y) \hspace{0.6cm} \mbox{as \,\,$n \to 0$},
\label{s3eq22}
\end{equation}
where $b$ is complex in order to match with the inner solution.
Thus, at $O(\frac 1 n)$ in (\ref{s3eq21}), we obtain
\begin{equation}
\tex{
      10 |{\mathrm e}^b| \left( \frac{{\mathrm d}b}{{\mathrm d}Y} \right)^9 + Y = 0  ,
      }
\label{s3eq23}
\end{equation}
whilst, at $O(1)$, we have
\begin{equation}
 \tex{
\big( 10 (b')^9 |{\mathrm e}^b| + \frac{Y}{10}\big) \frac{{\mathrm
d}B}{{\mathrm d}Y} + B \big( \alpha - \alpha \frac{Y}{10} b' +
(b')^8 |{\mathrm e}^b|\,
 (45 b'' + b'^2(1+\ln |B|) ) \big)
 }
=0, \label{s3eq24}
\end{equation}
where $'$ denotes $\frac {\mathrm d}{{\mathrm d}Y}$ and the approximation
\begin{equation}
 |f|^n \sim |{\mathrm e}^b| \left( 1 + n \ln |B|\right)
\label{s3eq25}
\end{equation}
has been used. The solution to (\ref{s3eq23}) that matches with (\ref{s3eq20}) is
\begin{equation}
 \tex{
  b(Y) = (1\pm i \lambda) \, 9 \, \ln \big( 1 - \cos \left(\frac{4\pi}{9} \right)
  \left(\frac{Y}{10} \right)^{\frac {10}9} \big),
  }
\label{s3eq26}
\end{equation}
with $\lambda$ as given in (\ref{s3lambda}). We can, in principle,
determine the amplitude $B(Y)$ via (\ref{s3eq24}). The finite
interface $y=y_0$, where $f$ vanishes, is thus given by
\begin{equation}
 \tex{
     y_0 \sim 10 \left( \frac{1}{n}\sec\left( \frac{4\pi}{9} \right) \right)^{\frac 9{10}}  \hspace{0.6cm}
     \mbox{ as\,\, $n \to 0$ } ,
     }
\label{s3eq27}
\end{equation}
again illustrating its divergent behaviour for small $n$.
\ssk

\setcounter{equation}{0}
\section{ Oscillatory solutions of changing sign near interfaces
via periodic structures}
\label{S4}


\noindent Here, we examine the local behaviour near the finite
interface $y_0=x_0/t^{\beta}$ for the similarity ODE  (\ref{ODE1})
in one dimension. We consider similarity profiles $f(y)$
exhibiting {\em maximal regularity} at the interface $y=y_0$, so
that, being extended by $f=0$ for $y > y_0$, these will give
solutions of the CP.

\subsection{Periodic structure of oscillations near interfaces for $n \in (0, \frac 98)$}

 For the thin film ODE (\ref{ODE1}), we have the asymptotic behaviour
 \begin{equation}
 \label{s4eq1}
 |f|^n f^{(9)} \sim \lambda_0 f \, , \quad \lambda_0 = \b y_0>0, \hspace{0.6cm} \mbox{as \,\, $y \to y_0^{-}$},
 \end{equation}
 where the no-flux condition in (\ref{i3}) has been used. To allow for
  an oscillatory behaviour, we seek solutions in the form
 \be
 \label{s4eq2}
 \mbox{$
 f(y) = (y_0-y)^{\mu} \phi(\eta) \,, \quad \eta =
 \ln(y_0-y), \quad \mbox{with} \,\,\, \mu=\frac{9}{n},
  $}
 \ee
 where the oscillatory component  $\varphi$
 satisfies the ninth-order autonomous ODE
 \be
 \label{eqLC}
 \tex{
 \sum_{k=0}^{9}  a_k \phi^{(9-k)} + \lambda_0|\phi|^{-n} \phi=0.
 }
 \ee
The coefficients $\{a_k\}$ are polynomials in $\mu$ of degree $k$,
namely
\[
  a_0 =1 , \hspace{0.5cm} a_1 = 9(\mu-4), \hspace{0.5cm} a_9 = \Pi_{i=0}^{8} (\mu-i),
\]
with the others easily obtainable using e.g., {\tt Maple} and not
recorded for conciseness.
\begin{center}
\begin{figure}[htp]
 \hspace{-1.5cm}
\includegraphics[scale=0.8]{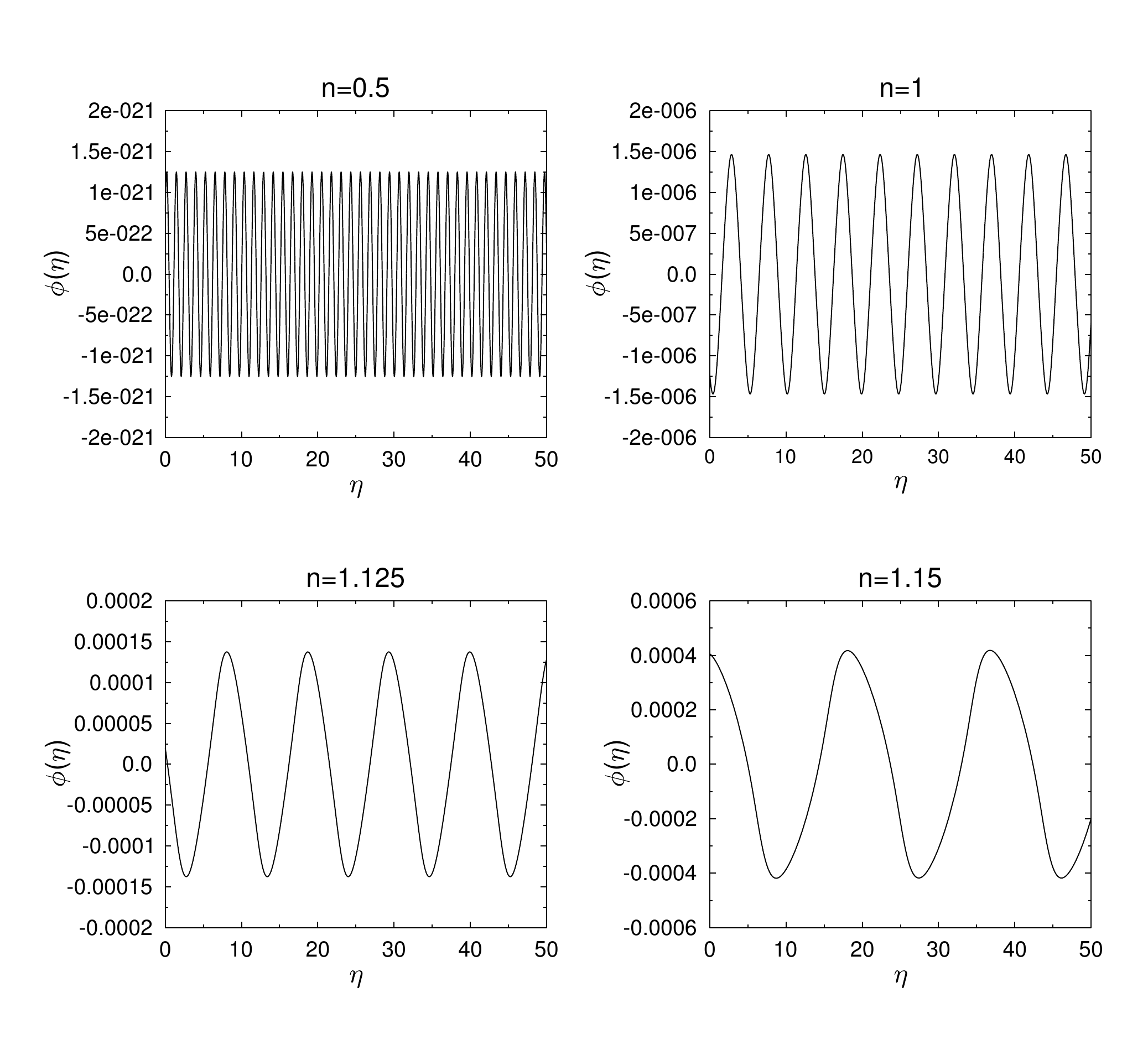}
\vskip -.4cm
 \caption{ \small Numerical illustration of the limit
cycles for selected $n$. In each case (\ref{eqLC}) with parameter
value $\lambda_0=1$ was solved as an IVP using {\tt MatLab} solver
{\tt ode15s}. Small error tolerances RelTol and AbsTol of
typically $10^{-13}$ were set, although these were relaxed for the
larger $n$ values.} \label{Fig1LC}
\end{figure}
\end{center}


\ssk

We formulate our overall (formal) understanding of the ODE
\ef{eqLC} as follows:

\ssk

\noi{\bf Conjecture \ref{S4}.1.}\,~{\em   For $n \in (0, \frac
98)$, the ODE $(\ref{eqLC})$ has a unique non-trivial
sign-changing periodic solution $\phi_*(\eta)$.
 }

\ssk

Thus, numerics suggest that this limit cycle is globally stable
and is unique (up to translations in $\eta$). Figure \ref{Fig1LC}
describes this stable periodic behaviour for selected  $n \in
(0,\frac 98)$. For $n \in (\frac{9}{8}, \frac{9}{7})$, global
stability fails since
 (\ref{eqLC})  admits also two equilibria $\phi=\pm \phi_0 $, where
 \be
 \label{Var55}
 \mbox{$
  \phi_0 =  [ -\frac{\lambda_0}{ \Pi_{i=0}^{8}(\mu-i)}]^{\frac 1n} >
  0.
 $}
 \ee
The amplitude of the periodic solution
decreases markedly as $n$ decreases, which suggests the need to
rescale for small $n$ as discussed later.

\subsection{Heteroclinic bifurcation of periodic solutions}

\begin{figure}[htp]
\centering
 \hskip -1cm
\includegraphics[scale=0.8]{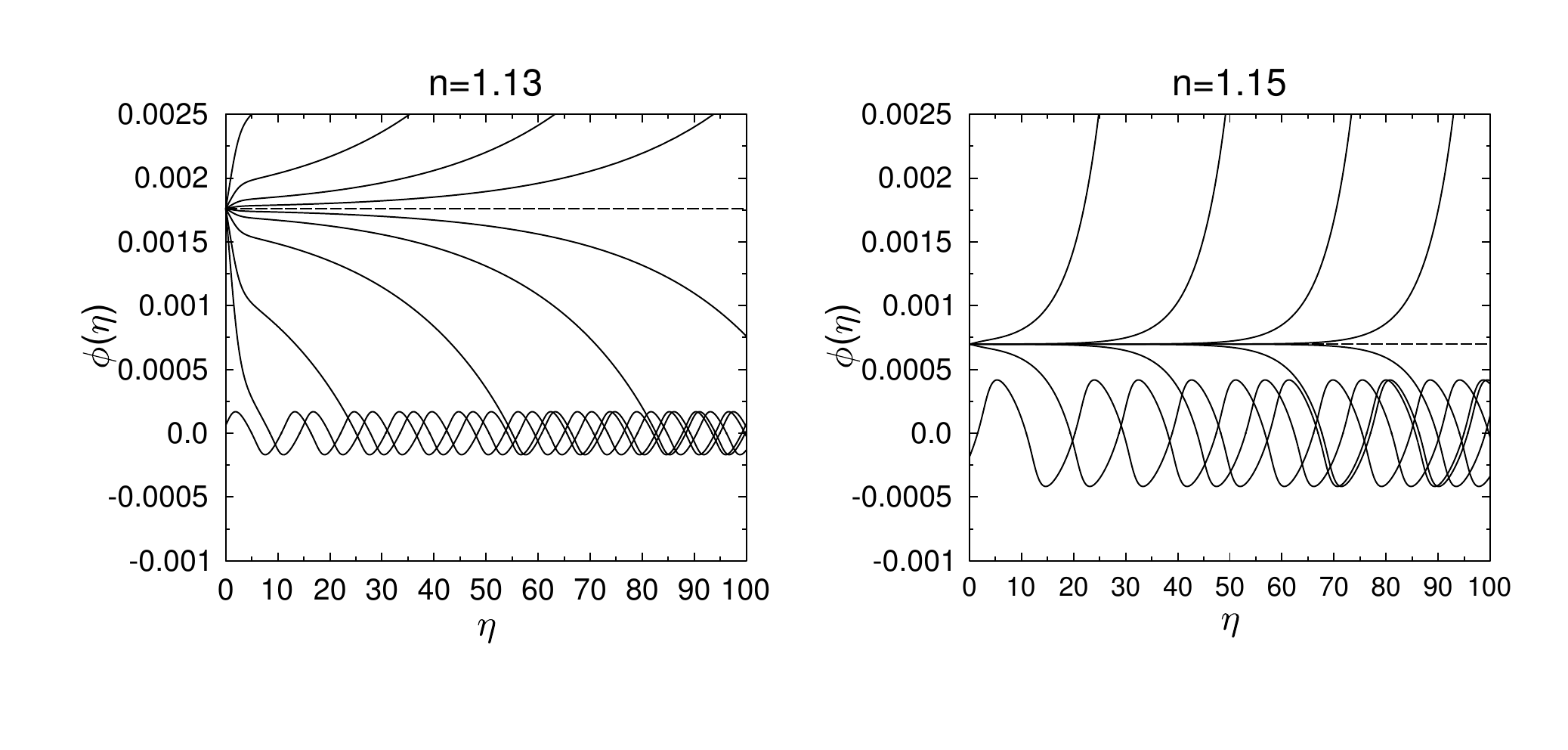}

\vskip -1cm \caption{ \small Illustration of numerical  solutions
to (\ref{eqLC}) with $\lambda_0=1$ for two selected values of
$n\in(\frac{9}{8},n_{\rm h})$. Shown are solutions leaving the
unstable constant solution $\phi=\phi_0$ with the globally stable
limit cycle. The figures are symmetric in $\phi$, with the
negative unstable constant solution $\phi=-\phi_0$ being omitted.}
\label{Fig2LC}
\end{figure}

It is crucial for both ODE and PDE theory to find a  precise
$n$-interval of existence of periodic, oscillatory solutions of
(\ref{eqLC}). Firstly, the stable periodic solution $\phi(\eta)$
persists to exist for $n> \frac 98$, where the constant solutions
$\phi=\pm \phi_0$ are unstable; see Figure \ref{Fig2LC} for
$n=1.13$ and $n=1.15$.


\begin{figure}[htp]
\centering
 \hskip -1cm
\includegraphics[scale=0.8]{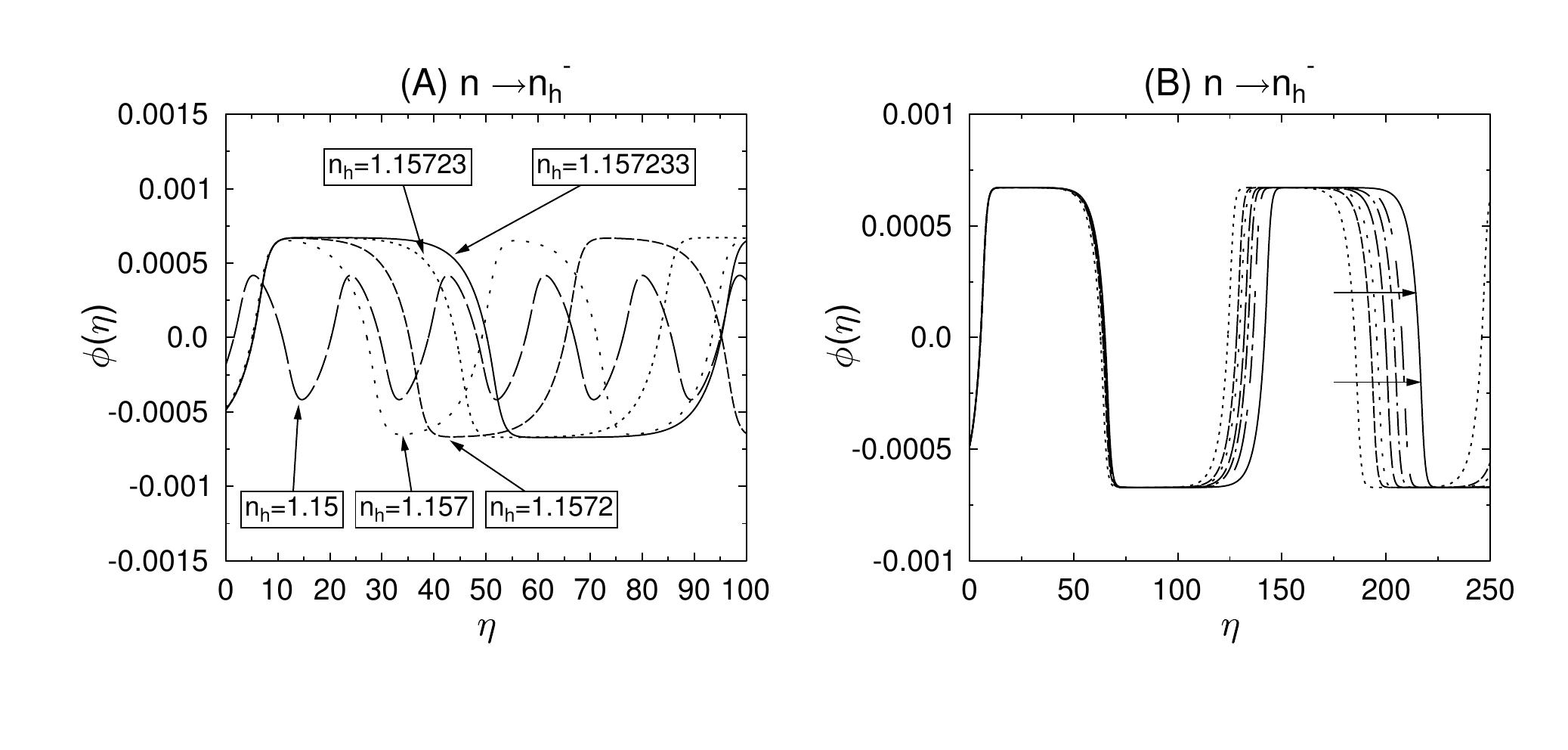}
\vskip -1cm \caption{ \small Formation of a heteroclinic
connection  $\phi_0 \to -\phi_0$ for the ODE (\ref{eqLC}),
$\lambda_0=1$, as $n\to n_{\rm h}^-$.
 }
\label{Fig3LC}
\end{figure}

Secondly, as $n$ increases further, the periodic solution
 is destroyed in a heteroclinic bifurcation, a phenomenon earlier observed for fourth- and sixth-order TFEs
  \cite{EGK2, EGK4}. The
following conjecture is entirely based on the numerical evidence.

\smallskip

\noi{\bf Conjecture \ref{S4}.2.}\,~{\em The stable periodic
solution of $(\ref{eqLC})$ exists for all $n \in (0,n_h)$, where
$n_h \in (\frac{9}{8}, \frac{9}{7})$ is a subcritical heteroclinic
$(\phi_0 \mapsto - \phi_0)$ bifurcation point of  stable periodic
solutions, which cease to exist for all $n \ge n_h$}.

\smallskip

Numerical calculations give
 \be
 \label{n**1}
 \fbox{$
 n_{\rm h}=1.1572339... \quad (\mbox{recall that $\frac 98=1.125$ and
 $\frac 97=1.2857...$}).
  $}
 \ee
Figure \ref{Fig3LC} shows formation of such a bifurcation as $n \to
n_h^-$. To obtain the bold line in Figure \ref{Fig3LC}(B), we took
$n=1.157233919$ (not all decimals being correct). This is a standard
scenario for homoclinic/heteroclinic
bifurcations, \cite[Ch.~4]{Perko}. A rigorous justification of such
non-local bifurcations is an still an open problem.

Thus, for
 $n$  larger
than $\frac 98$, not all the solutions are oscillatory near the
interfaces. For $n \in (\frac 98, \frac 97)$, there exists a
one-parametric bundle  of positive solutions with constant
$\phi(\eta)$ given by (\ref{Var55}). Nevertheless, for matching
purposes, the whole 2D asymptotic bundle
(\ref{s4eq2}) of oscillatory solutions has to be taken into
account, so that the oscillatory behaviour remains generic
(as in the linear case $n=0$ described next).

\subsection{Periodic solutions for small $n > 0$}

\begin{figure}[htp]
\centering
 \hskip -1.5cm
\includegraphics[scale=0.8]{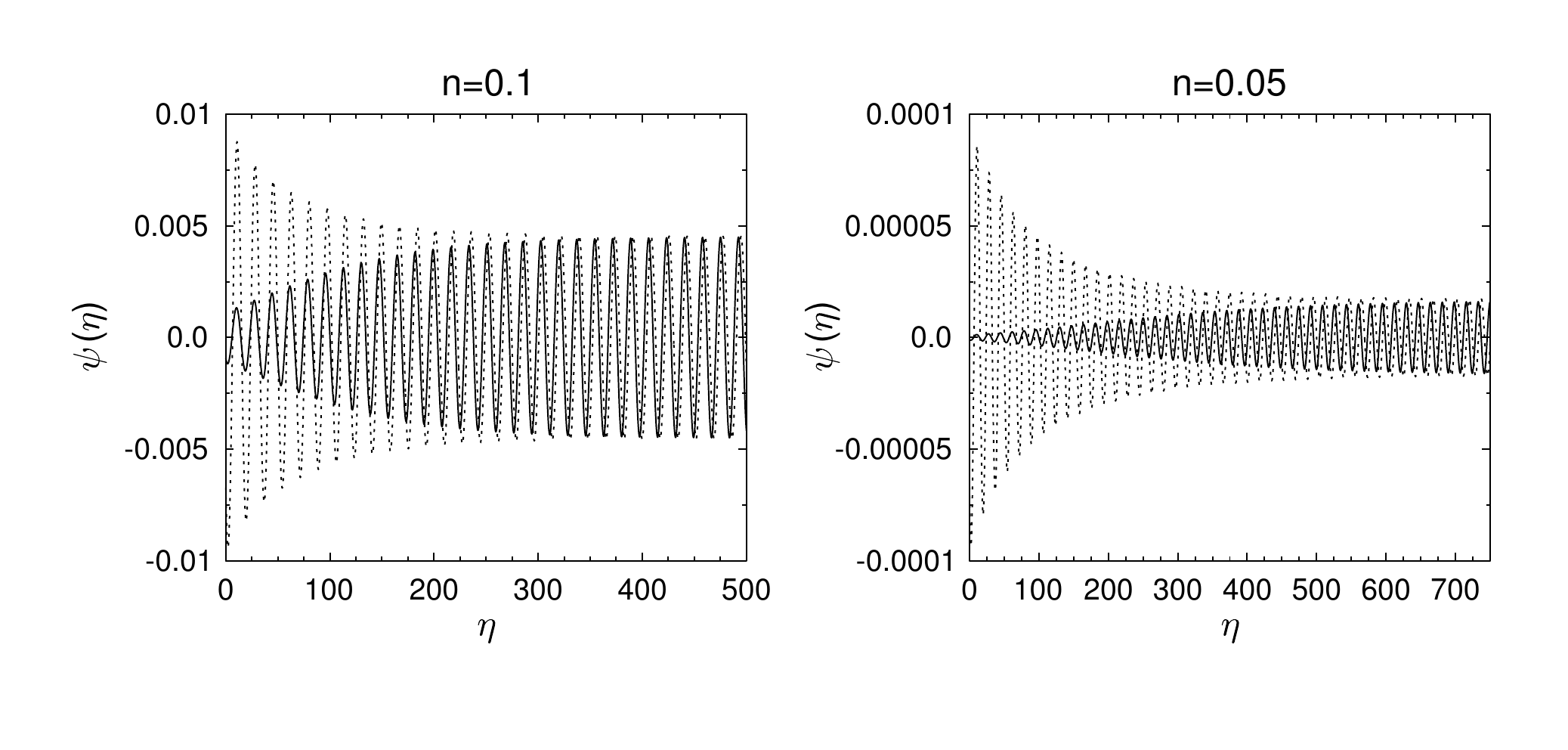}

\vskip -1cm \caption{ \small Stable periodic behaviour for the ODE (\ref{eqLC}), $\lambda_0=1$, for $n=0.1$ and $n=0.05$.
}
\label{Fig4LC}
\end{figure}

Here, we study the behaviour of periodic solutions for small
$n>0$. They are already difficult to detect by direct numerical
methods for $n=\frac 12$ as indicated in Figure \ref{Fig1LC},
where $\phi=O(10^{-21})$. To reveal the limiting oscillatory
behaviour as $n \to 0$, solutions of changing sign are expected to
be of the order
 $$
 \mbox{$
 \|\varphi\|_\infty \sim  (\frac n{9})^{\frac{9}n} \quad
 \mbox{for small} \,\,\, n>0.
  $}
  $$
 We thus rescale in equation (\ref{eqLC}) as follows
 \be
 \label{AAA.1}
 \mbox{$
\phi(\eta)= \left(\frac{n}{9} \right)^{\frac{9}{n}} \psi(s), \quad s = \frac{9}{n}  \eta .
 $}
 \ee
For small $n >0$, function $\psi(s)$ solves a simpler ODE with Euler's differential operator
and binomial coefficients, which can be
written in the form (omitting higher-order
 perturbations in $n$)
 \be
 \label{AAA.2}
\mbox{$
 {\mathrm e}^{-s} \frac{{\mathrm d}^9({\mathrm e}^s \psi)}{{\mathrm d}s^9}
 \equiv \sum_{k=0}^{9} \binom{9}{k} \psi^{(9-k)}
 =- \lambda_0 \frac{\psi}{|\psi|^n}.
$}
  \ee
Numerical analysis shows existence of a stable periodic solution
in (\ref{AAA.2}); see Figure \ref{Fig4LC}. It is worth mentioning that
the periodic oscillations become very small as $n$ reduces, e.g. by
 (\ref{AAA.1}),
 $$
 \|\phi \|_\infty \sim 10^{-176} \quad  \mbox{for $n=0.1$} .
  $$
Note that, for $n=\frac 12$, (\ref{AAA.1}) suggests $\|\phi
\|_\infty \sim 10^{-23}$ which is slightly smaller than the
numerical size indicated in Figure \ref{Fig1LC}. Stabilization to
periodic orbits of (\ref{AAA.2}) are shown in Figure \ref{Fig4LC}
for two sufficiently small $n$ values.

For $n=0$, the equation (\ref{AAA.2}) becomes linear,
\be
 \label{AAA.3}
\mbox{$
 \sum_{k=0}^{9} \binom{9}{k} \psi_0^{(9-k)}
 =- \lambda_0 \psi_0 .
$}
  \ee
with the characteristic equation $(\nu+1)^9+\lambda_0=0$ for
exponential solutions  $\psi(s)={\mathrm e}^{\nu s}$. Thus, the
generic stable behaviour for (\ref{AAA.3}) is exponential and
oscillatory:
\be
 \tex{
  \psi_0(s) \sim {\mathrm e}^{\big( \lambda_0^{\frac{1}{9}} \cos (\frac{\pi}{9}) -1 \big)s}
  \cos \big[ \big( \lambda_0^{\frac{1}{9}} \sin\frac{\pi}{9}\big)s
   + \mbox{constant} \big] \hspace{0.6cm} \mbox{as \,\, $s \to -\infty$}.
   }
\label{psi0}
 \ee
  These small $n$ asymptotics given by the scaling
(\ref{AAA.1}) describe the actual branching of periodic solutions
of (\ref{eqLC}) from the exponential decaying linear patterns
(\ref{psi0})
  for $n=0$.


\section{Towards global behaviour of nonlinear eigenfunctions via analytic approaches}
 \label{S5}

\ssk

\ssk


\subsection{Regularized problem}


To obtain global information about the solutions of the \emph{nonlinear eigenvalue problem} \eqref{self1}, i.e.,
its nonlinear eigenfunctions or source-type solutions of the
degenerate elliptic equation \eqref{self1},  we consider a homotopic deformation to the
 linear
elliptic problem
\begin{equation}
\label{s4}
 \tex{
    {\bf B}F \equiv -\D_y^5 F + \frac{1}{10}\, y \cdot \nabla_y F  +\frac{N}{10}\, F=0
    \quad \hbox{in} \quad \ren,\quad \int_{\ren} F(y) \, {\mathrm
    d}y=1.
    }
\end{equation}
This is the rescaled equation of the poly--harmonic equation  of
the tenth order \eqref{lin5}.  Note that \eqref{lin5} admits a unique classic
solution given by the convolution Poisson-type integral of the
form $$
 \tex{
    u(x,t)=b(t)\, * \, u^0 \equiv t^{-\frac{N}{10}} \int\limits_{\ren} F((x-z)t^{-\frac{1}{10}})
     u^0(z)\, {\mathrm d}z, \quad \mbox{where}
    }
$$
$$
 \tex{
    b(x,t)=t^{-\frac{N}{10}} F(y), \quad y:=\frac{x}{t^{1/10}} \quad  (x\in \ren)
    }
$$ is the unique fundamental solution of the operator
$\frac{\p}{\p t} - \D^{5}$ so that $F$ is the rescaled
fundamental kernel. Moreover, $F$ solves the linear elliptic equation
\eqref{s4} and operator ${\bf B}$ admits a countable family of
eigenfunctions
 (cf.
\cite{EGKP} for further details).

Therefore, we take the regularized uniform elliptic equation
\begin{equation}
\label{regul1}
 \tex{
     \nabla \cdot \left[ \phi_\e(f) \n \D^4 f\right] +\frac{1-\a n}{10}\, y \cdot \nabla f +\a
    f=0,
    \quad f \in C_0(\ren)\, , \quad \mbox{with}
    }
\end{equation}
 \be
  \label{phi2}
  \phi_\e(f)=|\e^2+f^2|^\frac n2, \quad \e>0,
   \ee
so that the inverse operator is smooth and analytic. Thus, for any
$\e \in (0,1]$, the uniformly elliptic equation admits a unique
classic solution $f=f_\e(y)$, which is an analytic function in
both variables $y$ and $\e$. Indeed, we would like to see under
which conditions we can have that $$ f_\e(y) \to f(y) \asA \e \to
0^+,
   $$ for a given well-defined  analytic functional family  (a
curve or a path),
 \be
 \label{Fem}
 {\mathcal P_\phi}=
\{f_\e(y)\}_{\e \in (0,1]}.
 \ee
  Thus, to obtain relevant information about the nonlinear eigenfunctions of the
  problem \eqref{self1} we will apply standard degree theory \cite{Kras,KZ} and first will perform
a kind of ``{\em double}" limit as  $\e,\, n \to 0^+$, where
special restrictions on the two parameters will be required.
Basically, because, passing to the limit just when $\e$ goes to
zero, we find a very deep problem since the regularized PDE loses
its uniform ellipticity.

Most of the existing results for \emph{thin film equations}
deal with {\em non-negative solutions} with compact support of
various FBPs, which are often more physically relevant and use standard integral identities for $\{f_\e\}$. In this
context, we should point out that such approximations for {\em
non-negative} and non-changing sign solutions,  with various
non-analytic (and non-smooth) regularizations (for example, of the
form $|u|^n +\e$, which is not analytic for $n<2$) have been
widely used before in TFE--FBP theory as a key foundation (cf. \cite{BF1})
but assuming the parabolic problem and using energy methods.
  Moreover, apart from the limiting problem when $n$ approximates
$0^+$  in the one-dimensional case it is not possible to apply the standard
  energy methods to ascertain the limiting behaviour in a convincing manner.
  Hence, we will use the degree theory and a kind of ``{\em double}" limit to resolve this issue.

Furthermore, although it looks quite reasonable to perform such a
limit when $\e \to 0^+$, as mentioned above just passing to this
limit we face many difficult problems since it is not sufficiently
clear, using for example integral identities techniques, how to
identify the limit (existence or non-existence of such a limit)
or,  even if we have more than one limit (see \cite{PV3}).
For this matter, we present a discussion
about how to deal with this particular limiting problem.



\subsection{Homotopy via degree theory}


First, we will perform a ``homotopy" transformation via standard degree theory when the
double limit $\e,\, n \to 0$. In fact, we will use
the existence of the limit when the parameters $n$ and $\e$ go to
zero in a certain manner.

In order to apply standard degree theory, we will write the
regularized equation \eqref{regul1} in the form
 \be
 \label{pertubeq}
 \tex{
     ({\bf B}_n+a{\rm Id})f_\e \equiv \D^5 f_\e +\frac{1-\a n}{10}\, y \cdot \nabla f_\e +(\a+a)
    f_\e= \nabla \cdot (1-\phi_\e(f_\e)) \n \D^4 f_\e + a f_\e,}
    \ee
        where $a>0$ is a parameter to be chosen so that the inverse
  operator $({\bf B}_n+a{\rm Id})^{-1}$ (a resolvent value) is a compact one in a
    weighted space $L^2_\rho(\ren)$, with $\rho$ a certain weight that makes the embedding of $H^{10}_\rho(\ren)$ compact into
  $L^2_\rho(\ren)$.
    Moreover, the
  spectrum of
    \be
\label{open}
\tex{{\bf B}_n+a{\rm Id} \equiv \D^5 +\frac{1-\a n}{10}\, y \cdot \nabla +(\a +a){\rm Id}}
\ee
    is always discrete and, actually, thanks to the spectral theory developed in \cite{EGKP}  for these higher-order operators, 
    for the operators
    \be
\label{oplin}
\tex{{\bf B}+a{\rm Id} \equiv \D^5 +\frac{1}{10}\, y \cdot \nabla +(\a+a){\rm Id},}
\ee
whose spectrum is
$$\tex{
    \s({\bf B})=\big\{\l_k := -\frac{k}{10}\,,\,k=0,1,2,...\big\},
    }
        $$
we have that
 \be
 \label{sp33}
  \tex{
 \s({\bf B}_n)=\big\{- \frac{k(1-\a n)}{10}+\a,
 \,k=0,1,2,...\big\},
 }
  \ee
 so that any choice of $a>0$ such that $a \not \in \s({\bf B}_n)$ is suitable in \ef{pertubeq}.

 We intend to perform a homotopy transformation from the \eqref{regul1} to the \eqref{s4} translating the already known oscillatory properties
of the self-similar poly-harmonic parabolic equation \eqref{lin5}
into the thin film equation \eqref{i1}. Note that since the
eigenfunctions of the elliptic equation \eqref{s4} are, up to
constant multipliers, derivatives of the fundamental kernel $F$
(and, for the adjoint one ${\bf B}^*$ are generalized Hermite
polynomials with finite oscillatory properties), our purpose will
be to get such an oscillation characteristic into the solutions of
the non-linear eigenvalue equation \eqref{self1} (the self-similar
thin film equation). Thus, using the degree theory and the
existence of convergence,
we ascertain some existence and multiplicity results for the
non-linear eigenvalue problem \eqref{self1}.

Homotopy deformations are used in other fields in mathematics,
especially geometry, to put in correspondence certain properties
of several geometrical objects and the topological degree is the
only invariant which is conserved by homotopic deformations.
However, we will use it as a tool to analyze topological
invariants of those geometrical objects which can be put in
correspondence with the considered equation providing us with a
natural method for studying the invariant properties of the
integral equation \eqref{pertubeq}.

First of all, through the next proposition we prove  that the
linear elliptic operator on the left hand side of the equation
\eqref{pertubeq} denoted by \eqref{open} ``converges", in a
natural sense, to the operator \eqref{oplin} when the parameter
$n$ goes to zero. Thus, it turns out that, when the parameter $n$
approximates zero, we have according to \ef{alb1} that
 $$
  \tex{
 \a_0(0)=\frac{N}{10}.
 }
 $$
  Moreover, extending that
approximation also for any $k \geq 1$,
the parameter $\a$  reaches the following family of values:
\begin{equation}
\label{bf4}
 \tex{
     \a_k(0) := -\l_k + \frac{N}{10} \quad \hbox{for any} \quad k=1,2,\ldots,
     }
\end{equation}
where $\l_k$ are the eigenvalues of the operator ${\bf B}$,
so that
\begin{equation*}
 \tex{
    \a_0(0)= \frac{N}{10}, \; \a_1(0)= \frac{N+1}{10},\; \a_2(0)= \frac{N+2}{10},\ldots ,
    \a_k(0)= \frac{N+k}{10}\ldots \,.
     }
\end{equation*}
Then, we introduce the next expression for the parameter $\a$
\begin{equation}
\label{i52}
    \tex{ \a_k(n):= \frac{N}{10+Nn}-\l_k.}
\end{equation}

\begin{proposition}
The operators \eqref{open}
$$  \tex{{\bf B}_n+a{\rm Id} \equiv \D^5 +\frac{1-\a_k(n) n}{10}\, y \cdot \nabla +(\a_k(n) +a){\rm Id}}
  $$
  converge to the operator \eqref{oplin}
   $$
    \tex{{\bf B}+a{\rm Id} \equiv \D^5 +\frac{1}{10}\, y \cdot \nabla +(\frac{N+k}{10}
        +a){\rm Id},}
    $$
   as $n\rightarrow 0$, in the generalized sense of Kato.
\end{proposition}

 \noi{\em Proof.} Indeed, for each $u\in H_0^{10}(B_1)$ we have that
 $$
 \tex{ \left\| ({\bf B}_n+a{\rm Id}) u- ({\bf B}+a{\rm Id})u\right\|_{L^2(B_1)}
 \leq n \left\| \a_k(n) y\nabla u \right\|_{L^2(B_1)}.}
 $$
 Hence, from the expression for the parameter $\a_k(n)$ and Sobolev's inequality
 $$
 \tex{ \left\| ({\bf B}_n+a{\rm Id}) u- ({\bf B}+a{\rm Id})u\right\|_{L^2(B_1)}
 \leq c K \left\| u\right\|_{H_0^{10}(B_1)},}
 $$
 with $K>0$, a positive constant.
 Therefore, for any $\e >0$, there exists $n_0$ such that
 $$
 \tex{ \left\| ({\bf B}_n+a{\rm Id}) u- ({\bf B}+a{\rm Id})u\right\|_{L^2(B_1)}
 \leq \e \left\| u\right\|_{H_0^{10}(B_1)},}
 $$
 for all $n\in (0,n_0)$ and  $u\in H_0^{10}(B_1)$.
 $\qed$

\vspace{0.2cm}

Subsequently, using the compact embedding of $H^{10}_\rho(\ren)$ into
  $L^2_\rho(\ren)$,  we find that.
  \begin{proposition}
\be
\label{converep}
\tex{f_\e \longrightarrow \hat{F},}
\ee
performing a double limit as $n$ and $\e$ go to zero, at least in $L^2_\rho(\ren)$.
\end{proposition}
 \noi{\em Proof.}
 So far, we cannot identify which problem $\hat{F}$ belongs to.
However, we write the equation \eqref{pertubeq} in the integral form
\be
 \label{pertubeq11}
 \tex{
    f_\e= ({\bf B}_n+a{\rm Id})^{-1} \left[ \nabla \cdot (1-|\e^2+f_\e^2|^\frac n2) \n \D^4 f_\e +a f_\e \right],}
    \ee
with ${\bf B}_n+a\rm{Id}$ denoted by \eqref{open},  for which we
know the expression for the whole spectrum explicitly and, also,
that this operator is compact in a weighted space $L^2_\rho(\ren)$
with the existence of the inverse for a suitable and positive
$a\notin \sigma({\bf B}_n)$. For the nonlinear term $$\nabla \cdot
(1-|\e^2+f_\e^2|^\frac n2 ) \n \D^4 f_\e +a f_\e,$$ we have that
it is relatively compact, thanks to the existence of convergence
shown previously \eqref{converep} and assuming the condition
\be
  \label{nn12}
   \tex{
    \mbox{for} \quad \d \sim \e, \quad
  n=n(\e) \to 0 \,\,\, \mbox{such that} \quad \e^{\frac{n(\e)}2}
  \to 0.
   }
   \ee
   Indeed, setting
\begin{equation*}
    \tex{\lim_{\e \rightarrow 0} \e^{\frac{n(\e)}2} = 0, \quad \hbox{then} \quad
     \lim_{\e \rightarrow 0} n(\e) \ln \e= -\iy.}
\end{equation*}
Hence, taking
\be
\label{expan}
\tex{F_\e(f_\e)= 1-|\e^2+f_\e^2|^\frac n2 =  - \frac n2 \, \ln (\e^2+f_\e^2)(1+o(1)) \asA n \to 0^+,}
\ee
for a family $\{f_\e(y)\}$ of uniformly bounded and smooth
 solutions, when $f_\e \approx 0$ yields the demand
  \be
  \label{keycond}
   \fbox{$
  n \, |\ln \e(n)| \to 0 \asA n \to 0,
   $}
  \ee
 which it is true if we assume \ef{nn12} such that the regularization parameter $\e \ll {\mathrm e}^{- \frac 1n}$.
Thus, substituting \eqref{expan} into \eqref{pertubeq11}, we
arrive at
\be
\label{pertubeq12}
\tex{
    f_\e= ({\bf B}_n+a{\rm Id})^{-1} \left[ \nabla \cdot (- \frac n2 \, \ln (\e^2+f_\e^2)(1+o(1))) \n \D^4 f_\e +a f_\e \right],}
    \ee
    and passing to the limit when $n$ and $\e(n)$ go to zero we find that
    there exists a fixed point for the integral equation
   \be
 \label{fixedpoeq}
 \tex{
    \hat{F}= ({\bf B}+a{\rm Id})^{-1} (a \hat{F}),}
    \ee
    whose solutions are the eigenfunctions of the linear elliptic problem \eqref{oplin}, i.e., $\hat{F}=\psi$.
    $\qed$

\vspace{0.2cm}

        Note that \eqref{expan} could be replaced by a more general form such as \eqref{s3eq22} with
        $$\tex{G_\e(f_\e)= \frac{b}{n}+\ln B+ O(n),\asA n \to 0^+,}$$
        obtaining a different condition from \eqref{keycond}.

        Therefore, applying the degree theory, together with Fixed Point
        Theory
    (see \cite{Kras,KZ} for any further details), we can assure the existence of a countable family of a  direct $n$-expansion of the
   solutions for the problem \eqref{self1} to guarantee branching at $n=0^+$. In fact, the degree provides us with the existence
   of continuous branches of eigenfunctions for the equation \eqref{self1} since it stays invariant via homotopic deformations as the ones performed here. In
   the terms exposed by Krasnosel'skii \cite{Kras}, we would talk about the rotation of the vector field  of the form
   $$\Phi=\rm{Id} - {\bf G},$$
   where ${\bf G}$ is the operator on the right hand side of \eqref{pertubeq12}. Note that
    the invariance analysis of the rotation of vector fields and
   the degree in the sense of Leray--Schauder are equivalent.


\subsection{The limiting problem just when $\e \to 0^+$}


In general, 
to ascertain the limit $$ \tex{f_\e
\longrightarrow \hat{F},} $$
 just when $\e \to 0^+$ instead of the double limit performed
above, we take into account that the inverse operator
 $$\tex{({\bf
B}_n+a{\rm Id})^{-1}}
  $$
   is smooth and analytic for any $n$, and
the convergence of the sequence $\{f_\e\}$, at least, in
$L^2_\rho(\ren)$. Here, again, $a>0$ is a parameter to be chosen
so that the inverse
  operator is compact in the weighted space $L^2_\rho(\ren)$, with $\rho$ being  a certain weight that
   makes the embedding of $H^{10}_\rho(\ren)$ compact into
  $L^2_\rho(\ren)$. Moreover, the solutions of the regularized problem \eqref{regul1},
\eqref{phi2} are analytic in both variables $\e$, $y$, and, by the
construction, we know that $$\tex{f_\e \in C_0(\re^N),}$$
 i.e., these
solutions have compact support and, also, by the conservation of
mass, $$\tex{\int_{\re^N} |f_\e|\leq C, \quad \hbox{with}\quad
C>0, \,\,\, \hbox{being a positive constant}.}
  $$

Observe that we have solutions of changing sign, with compact
support and exponential decay. Then, together with the boundary
conditions, we find that 
  $$
  \tex{\int_{\re^N} |f_\e|^2 =
\int_{\re^N\setminus \{|f_\e|<\d \}} |f_\e|^2 + \int_{\{|f_\e|<\d
\}} |f_\e|^2 \leq M+ \d^2 |\{\supp f_\e\}|,} 
  $$ 
  providing us with an
estimate for the norms in $L^2$. Here, we assume that the
solutions are oscillatory of changing sign, but we are not able to
extract an information about those oscillations (since the
argument we have done before to extract information from the
solutions at $n=0$ is not applicable now), since we do not posses
an {\it a priori} information about this oscillatory property when
$n\neq 0$. However, the goal would be to extend analytically these
oscillatory properties from $n=0$ forward.

Furthermore, since the inverse operator of \eqref{open} is smooth
and analytic, we can assure that \eqref{open} is a topological
isomorphism and, hence, we can apply the Implicit Function Theorem
to the equation \eqref{pertubeq} with the parameter $n$ fixed.

Therefore, it looks quite natural to apply the argument of passing
to the limit just as $\e$ goes to zero in the equation
\eqref{pertubeq11}, $$
 \tex{
    f_\e= ({\bf B}_n+a{\rm Id})^{-1} \left[ \nabla \cdot (1-|\e^2+f_\e^2|^\frac n2) \n \D^4 f_\e +a f_\e \right].}
    $$
 However, we face here several
problems that make this final process very tricky. Indeed, to get
the convergence of the previous perturbed equation
\eqref{pertubeq}, we need to get the term
    \be
    \label{expterm}
    \tex{\nabla \cdot (1-|\e^2+f_\e^2|^\frac n2) \n \D^4 f_\e + af_\e,}
    \ee
    bounded in $L^2_\rho(\ren)$. Essentially, since we have an inverse compact operator $({\bf B}_n+a{\rm
    Id})^{-1}$,
    if \eqref{expterm} is bounded in $L^2_\rho(\ren)$, we might be able to find a convergent subsequence being the solutions $f_\e$
    of that fixed point equation relatively compact. Nevertheless, on the contrary from what we had above
     for the double limit, here, we cannot assure that
    the non-linear term \eqref{expterm} is relatively compact making the solution of the Fixed Point equation \eqref{pertubeq11} far from obvious.
    Even though, it looks quite reasonable.
    
    Moreover, when $|f_\e|\geq \d$ for $\d>0$, since the solutions of the perturbed equation \eqref{pertubeq}
    are continuous with compact support, we have that the $L^\infty$-norm is bounded in those subsets
    $$\tex{\|f_\e\|_{L^\infty,\{|f_\e|\geq \d\}} < K, \quad \hbox{for a constant}\quad K>0,}$$
    and, hence, the convergence of the fixed point equation \eqref{pertubeq} when $\e \to 0^+$ is
    guaranteed, at least, for the particular $n$'s, for which there exists a solution obtained numerically in other sections.

    However, that is not so clear to obtain, at least directly, when $|f_\e|< \d$ for $\d>0$.
    Indeed, by construction,
    we find that
    $$\tex{ ({\bf B}_n+a{\rm Id} )f_\e\in L_\rho^2 (\re^N),}$$
    with $f_\e\in H_\rho^1(\re^N)$, but we cannot imply directly that
    $$\tex{\nabla \cdot (1-|\e^2+f_\e^2|^\frac n2) \n \D^4 f_\e +af_\e \in L_\rho^2 (\re^N),} $$
    from the perturbed equation \eqref{pertubeq}.
    Moreover, computing
    $$\tex{ \int_{\re^N} (\phi_\e(f_\e))^2 (\n \D^4 f_\e)^2= \int_{\re^N} |\e^2+f_\e^2|^n (\n \D^4 f_\e)^2
    = \int_{\re^N} |\e^2+f_\e^2|^\frac n2 |\e^2+f_\e^2|^\frac n2 (\n \D^4 f_\e)^2,}$$
    we arrive at
    $$\tex{ \int_{\re^N} (\phi_\e(f_\e))^2 (\n \D^4 f_\e)^2 \leq K,}$$
    for a positive constant $K$, assuming that
    \be
    \label{firest}
    \tex{\int_{\re^N}  |\e^2+f_\e^2|^\frac n2 (\n \D^4 f_\e)^2\leq K \quad \hbox{and}\quad \int_{\re^N} |f_\e|\leq K.}
    \ee
    This final argument would provide us with the convergence of the fixed point equation, but, to this aim, we need
    the first estimation of \eqref{firest} and the relatively compactness of the non-linear terms \eqref{expterm}.

    Finally, even though we ascertain the existence of the limit of the equation \eqref{pertubeq11}, we cannot assure
    how many solutions satisfy the limiting problem
    $$
 \tex{
    \hat{F}= ({\bf B}_n+a{\rm Id})^{-1} \left[ \nabla \cdot (1-|\hat{F}|^n) \n \D^4 \hat{F} +a \hat{F} \right].}
    $$
 and also, what kind of solutions are, since we do not have an {\it a priori}
information about the solutions for other $n$ apart from $n=0$ (as performed above).

\section{Nonlinear eigenfunctions: numerical approach}
\label{S6}

\begin{figure}[htp]

\hspace*{-1cm}
\includegraphics[scale=0.8]{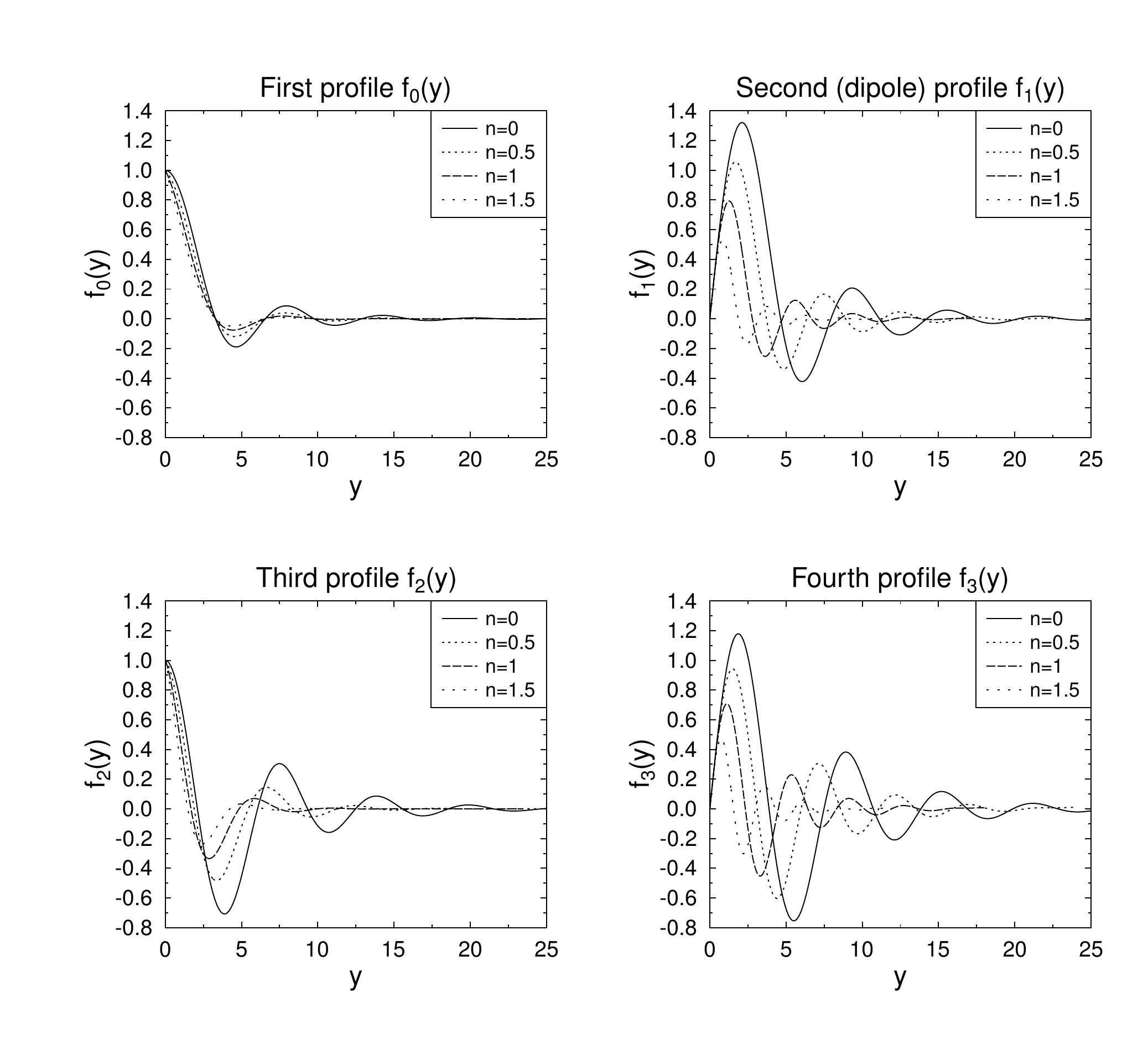}

\vskip -0.5cm \caption{ \small Numerical solution of the first four nonlinear eigenfunctions profiles $f_k,k=0,1,2,3$ for selected $n$ in one-dimension $N=1$.  }
 \label{ef}
\end{figure}

\noindent Finally, we construct numerically the nonlinear eigenfunctions in one space dimension. The nonlinear eigenvalue problem (\ref{self1}) for $N=1$ becomes
 \begin{equation}
\label{self1N1}
 \tex{
     \left( |f|^{n} f^{(9)}\right)' +\frac{1-\a n}{10}\, y f' +\a
    f=0,
    \quad f \in C_0(\re)\, ,
    }
\end{equation}
 for $n>0$, with $f$ being, thus,  compactly supported. For $n=0$, instead, we naturally  require $f$ to have exponential decay in infinity, now belonging to an appropriately weighted
 $L^2$-space as stated in (\ref{WW11}). The nonlinear eigenvalue-eigenfunction pairs are denoted by
  $\{\alpha_k(n),f_k\}$ for $k=0,1,2,3,\, ...$, and the eigenfunctions
 are normalised using
 \begin{equation}
 f_k(0)=1, \hspace{0.25cm} k=0,2,4,\ldots;\,\,\, f_k'(0)=1, \hspace{0.5cm} k=1,3,5,\ldots.
 \end{equation}
The first eigenvalue-eigenfunction pair $\{\alpha_0(n),f_0\}$
preserves mass, so that (\ref{self1N1}) may be integrated once to
give
\begin{equation}
\label{f0}
 \tex{
     |f_0|^{n} f_0^{(9)}   +  \a_0 y f_0 =0, \hspace{0.5cm} \mbox{with} \;\;\; \alpha_0(n)=\frac{1}{10+n} \;,
    }
\end{equation}
and is completed with the boundary conditions
\begin{eqnarray}
 \mbox{at $y=0$:} \hspace{0.6cm} &&  f_0=1,  f_0^{(i)}=0 \;\;\;\mbox{for $i=1,3,5,7$}, \label{f0bc1} \\
 \mbox{at $y=y_0$:} \hspace{0.6cm} &&  f_0= f_0^{(i)} = 0 \;\;\;\mbox{for $i=1,2,3,4$}. \label{f0bc2}
\end{eqnarray}
Since $\alpha_0$ is known, this gives a tenth-order system when $n>0$ to determine $f_0$ and the finite free boundary $y_0>0$ (the corresponding interface
being $x=y_0 t^{\beta_0}$ with $\beta_0$ as given in \eqref{sf3}). When $n=0$, then $y_0=\infty$.
Figure \ref{ef} shows illustrative $f_0$ profiles for selected $n$ values in one-dimension (N=1).

The system was solved as an IVP in {\tt MatLab} (shooting from
$y=0$), using the ODE solver ode15s with error tolerances of
AbsTol=RelTol=$10^{-10}$ and the regularisation $|f|^n \mapsto
(f^2+\delta^2)^{n/2}$ with $\delta=10^{-10}$. Since the $\alpha_0$
 is explicitly  known, this gives a tenth-order system when $n>0$ to determine
$f_0$ and the finite free boundary $y_0$. When $n=0$, then
$y_0=\infty$.

\begin{figure}[htp]

\hspace*{-1cm}
\includegraphics[scale=0.4]{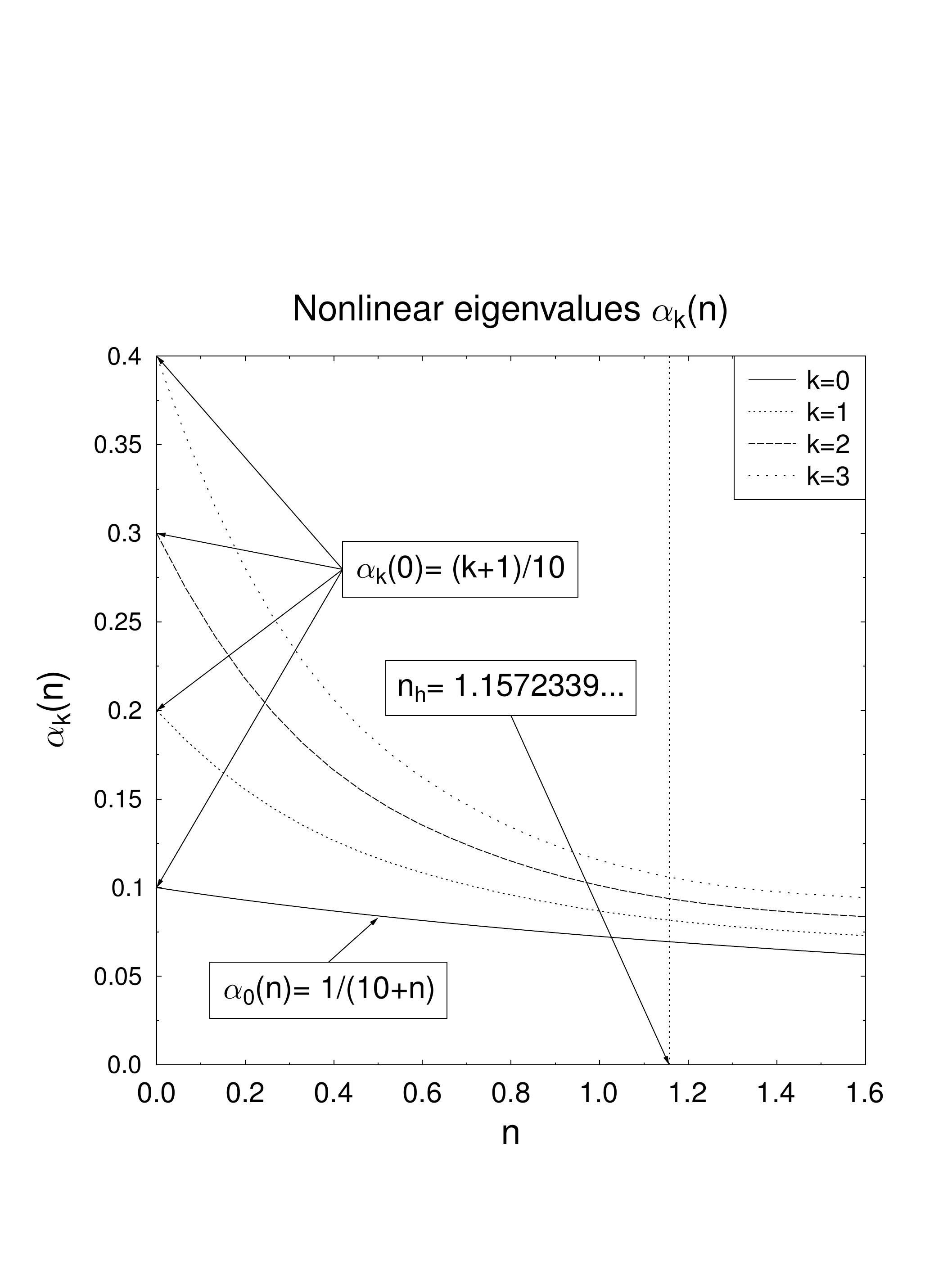}

\vskip -0.5cm \caption{ \small Numerical construction of the first four nonlinear eigenvalues $\alpha_k(n)$ in one-dimension $N=1$. }
 \label{ev}
\end{figure}

The other eigenvalue-eigenfunction pairs $\{\alpha_k,f_k\}$ for
$k\geq 1$  and $n>0$ satisfy the ODE in (\ref{self1N1}) with
\begin{eqnarray}
 \mbox{at $y=0$:} \hspace{0.5cm} &&  \left\{ \begin{array}{ll}
                             f_k=1, \hskip 0.25cm f_k^{(i)} =0 \;\;\;\mbox{for $i=1,3,5,7,9$},
                             \hspace{0.5cm} & \mbox{if $k$ is even,}\label{fkbc1a} \\
                             f_k'=1, \hskip 0.25cm f_k=f_k^{(i)}=0 \;\;\;\mbox{for $i=2,4,6,8$}, \hspace{0.5cm} &\mbox{if $k$ is odd} \label{fkbc1b} \\
                             \end{array} \right.
\end{eqnarray}
and
\begin{eqnarray}
\mbox{at $y=y_0$:} \hspace{0.5cm} &&  f_k = f_k^{(i)} =0 \;\;\;\mbox{for $i=1,2,3,4,5$}. \label{fkbc2}
\end{eqnarray}

Figure \ref{ef} show the eigenfunction profiles for the first four cases $k=0,1,2,3$, obtained by using the same shooting numerical
procedure for the first profile (but appropriately adapted for this 12th-order system). A plot of the eigenvalues is given in Figure\;\ref{ev}.

The case when $n=0$ requires a
slight modification, with $y_0=\infty$ and is discussed in \cite{AEG1}.


\end{document}